\documentclass[12pt]{article}
\scrollmode
\begin{document}
\title{\bf Random walks with strongly inhomogeneous rates
and singular diffusions:\\ convergence, localization
and aging in one dimension}
\author{L.~R.~G.~Fontes \and M.~Isopi \and C.~M.~Newman}
%\date{}
\maketitle

\newtheorem{defin}{Definition}[section]
\newtheorem{Prop}{Proposition}
\newtheorem{teo}{Theorem}[section]
\newtheorem{ml}{Main Lemma}
\newtheorem{con}{Conjecture}
\newtheorem{cond}{Condition}
\newtheorem{prop}{Proposition}[section]
\newtheorem{lem}{Lemma}[section]
\newtheorem{rmk}{Remark}[section]
\newtheorem{cor}{Corollary}[section]
\renewcommand{\theequation}{\thesection .\arabic{equation}}

\newcommand{\beq}{\begin{equation}}
\newcommand{\eeq}{\end{equation}}
\newcommand{\beqn}{\begin{eqnarray}}
\newcommand{\beqnn}{\begin{eqnarray*}}
\newcommand{\eeqn}{\end{eqnarray}}
\newcommand{\eeqnn}{\end{eqnarray*}}
\newcommand{\bprop}{\begin{prop}}
\newcommand{\eprop}{\end{prop}}
\newcommand{\bteo}{\begin{teo}}
\newcommand{\bcor}{\begin{cor}}
\newcommand{\ecor}{\end{cor}}
\newcommand{\bcon}{\begin{con}}
\newcommand{\econ}{\end{con}}
\newcommand{\bcond}{\begin{cond}}
\newcommand{\econd}{\end{cond}}
\newcommand{\eteo}{\end{teo}}
\newcommand{\brm}{\begin{rmk}}
\newcommand{\erm}{\end{rmk}}
\newcommand{\blem}{\begin{lem}}
\newcommand{\elem}{\end{lem}}
\newcommand{\ben}{\begin{enumerate}}
\newcommand{\een}{\end{enumerate}}
\newcommand{\bei}{\begin{itemize}}
\newcommand{\eei}{\end{itemize}}
\newcommand{\bdf}{\begin{defin}}
\newcommand{\edf}{\end{defin}}

\newcommand{\nn}{\nonumber}
\renewcommand{\=}{&=&}
\renewcommand{\>}{&>&}
\renewcommand{\le}{&\leq&}
\newcommand{\+}{&+&}
\newcommand{\fr}{\frac}

\newcommand{\Bbb}{\bf}
\renewcommand{\r}{{\Bbb R}}
\newcommand{\Z}{{\Bbb Z}}
\newcommand{\z}{{\Bbb Z}}
\newcommand{\zd}{\z^d}
\newcommand{\zz}{{\Bbb Z}}
\newcommand{\R}{{\Bbb R}}
\newcommand{\E}{{\Bbb E}}
\newcommand{\C}{{\Bbb C}}
\renewcommand{\P}{{\Bbb P}}
\newcommand{\N}{{\Bbb N}}
\newcommand{\var}{{\Bbb V}}
\renewcommand{\S}{{\cal S}}
\newcommand{\T}{{\cal T}}
\newcommand{\cm}{{\cal M}}
\newcommand{\cp}{{\cal P}}
\newcommand{\xt}{X_t}
\renewcommand{\ge}{g^{(\epsilon)}}
\newcommand{\xe}{y^{(\epsilon)}}
\newcommand{\ye}{y^{(\epsilon)}}
\newcommand{\bx}{{\bar y}}
\newcommand{\by}{{\bar y}}
\newcommand{\bw}{{\bar w}}
\newcommand{\bxe}{{\bar y}^{(\epsilon)}}
\newcommand{\bye}{{\bar y}^{(\epsilon)}}
\newcommand{\bwe}{{\bar w}^{(\epsilon)}}
\newcommand{\bxz}{{\bar y}}
\newcommand{\bwz}{{\bar w}}
\newcommand{\we}{w^{(\epsilon)}}
\newcommand{\Xe}{Y^{(\epsilon)}}
\newcommand{\Ze}{Z^{(\epsilon)}}
\newcommand{\Ye}{Y^{(\epsilon)}}
\newcommand{\tye}{{\tilde Y}^{(\epsilon)}}
\newcommand{\hy}{{\hat Y}}
\newcommand{\ve}{V^{(\epsilon)}}
\newcommand{\Ne}{N^{(\epsilon)}}
\newcommand{\ce}{c^{(\epsilon)}}
\newcommand{\cle}{c^{(\l\epsilon)}}
\newcommand{\xet}{Y^{(\epsilon)}_t}
\newcommand{\hxt}{\hat X_t}
\newcommand{\btn}{\bar\tau_n}
\newcommand{\ct}{{\cal T}}
\newcommand{\rn}{{\cal R}_n}
\newcommand{\nt}{{N}_t}
\newcommand{\lnk}{{\cal L}_{n,k}}
\newcommand{\cl}{{\cal L}}
\newcommand{\tc}{\tilde{\cal C}_b}
\newcommand{\hxtt}{\hat X_{\ct}}
\newcommand{\txnt}{\tilde X_{\nt}}
\newcommand{\xs}{X_s}
\newcommand{\xn}{\tilde X_n}
\newcommand{\tx}{\tilde X}
\newcommand{\hx}{\hat X}
\newcommand{\txi}{\tilde X_i}
\newcommand{\txij}{\tilde X_{i_j}}
\newcommand{\taxi}{\tau_{\txi}}
\newcommand{\txn}{\tilde X_N}
\newcommand{\xk}{X_K}
\newcommand{\ts}{\tilde S}
\newcommand{\im}{I^-}
\newcommand{\ip}{I^+}
\newcommand{\hal}{H_\a}
\newcommand{\ba}{B_\a}

\renewcommand{\a}{\alpha}
\renewcommand{\b}{\beta}
\newcommand{\g}{\gamma}
\newcommand{\G}{\Gamma}
\renewcommand{\d}{\delta}
\newcommand{\D}{\Delta}
\newcommand{\e}{\epsilon}
\newcommand{\fes}{\phi^{(\epsilon)}_s}
\newcommand{\fet}{\phi^{(\epsilon)}_t}
\newcommand{\fe}{\phi^{(\epsilon)}}
\newcommand{\pset}{\psi^{(\epsilon)}_t}
\newcommand{\pse}{\psi^{(\epsilon)}}
\renewcommand{\l}{\lambda}
\newcommand{\me}{\mu^{(\epsilon)}}
\newcommand{\re}{\rho^{(\epsilon)}}
\newcommand{\tre}{\tilde{\rho}^{(\epsilon)}}
\newcommand{\nue}{\nu^{(\epsilon)}}
\newcommand{\mbe}{{\bar\mu}^{(\epsilon)}}
\newcommand{\rbe}{{\bar\rho}^{(\epsilon)}}
\newcommand{\mb}{{\bar\mu}}
\newcommand{\rb}{{\bar\rho}}
\newcommand{\mbz}{{\bar\mu}}
\newcommand{\s}{\sigma}
\renewcommand{\o}{\Omega}
\newcommand{\tio}{\tilde\S}
\renewcommand{\sl}{\sigma'}
\newcommand{\si}{\s(i)}
\newcommand{\sit}{\s_t(i)}
\newcommand{\ei}{\eta(i)}
\newcommand{\eit}{\eta_t(i)}
\newcommand{\eot}{\eta_t(0)}
\newcommand{\sil}{\s'_i}
\newcommand{\sj}{\s(j)}
\newcommand{\st}{\s_t}
\newcommand{\so}{\s_0}
\newcommand{\xii}{\xi_i}
\newcommand{\xij}{\xi_j}
\newcommand{\xio}{\xi_0}
\newcommand{\ti}{\tau_i}
\newcommand{\te}{\tau^{(\epsilon)}}
\newcommand{\bt}{\bar\tau}
\newcommand{\tti}{\tilde\tau_i}
\newcommand{\tto}{\tilde\tau_0}
\newcommand{\tei}{T_i}
\newcommand{\ttei}{\tilde T_i}
\newcommand{\tes}{T_S}
\newcommand{\tao}{\tau_0}
\renewcommand{\t}{\tilde t}

\newcommand{\da}{\downarrow}
\newcommand{\ua}{\uparrow}
\newcommand{\ar}{\rightarrow}
\newcommand{\va}{\stackrel{v}{\rightarrow}}
\newcommand{\ppa}{\stackrel{pp}{\rightarrow}}
\newcommand{\dw}{\stackrel{w}{\Rightarrow}}
\newcommand{\Va}{\stackrel{v}{\Rightarrow}}
\newcommand{\Ppa}{\stackrel{pp}{\Rightarrow}}
\newcommand{\lar}{\leftrightarrow}
\newcommand{\la}{\langle}
\newcommand{\ra}{\rangle}
\newcommand{\ep}{\vspace{.5cm}}
\newcommand\square{\ifmmode\sqr\else{$\sqr$}\fi}
\newcommand\sqr{\vcenter{
          \hrule height.1mm
          \hbox{\vrule width.1mm height2.2mm\kern2.18mm\vrule width.1mm}
          \hrule height.1mm}}                  % This is a slimmer sqr.

%%%%%%%%%%%%%%%%%% ABSTRACT %%%%%%%%%%%%%%%%%%%%%%%%%%%%%%%%
\begin{abstract}

Let $\tau = (\tau_i : i \in {\Bbb Z})$ denote i.i.d.~positive
random variables with common distribution $F$ and
(conditional on $\tau$) let
$X\, = \,(X_t : t\geq0,\,X_0=0)$, be a continuous-time simple
symmetric random walk on ${\Bbb Z}$ with inhomogeneous rates
$(\tau_i^{-1} : i \in {\Bbb Z})$. When $F$ is in the domain of attraction
of a stable law of exponent $\alpha<1$ (so that ${\Bbb E}(\tau_i) = \infty$
and X is subdiffusive), we prove that $(X,\tau)$, suitably
rescaled (in space and time), converges to a natural (singular)
diffusion $Z \, = \, (Z_t : t\geq0,\,Z_0=0)$ with a random (discrete) speed
measure
$\rho$. The convergence is such that the ``amount of localization'',
$\E \sum_{i \in {\Bbb Z}} [\P(X_t = i|\tau)]^2$ converges as $t \to \infty$ to
$\E \sum_{z \in {\Bbb R}} [\P(Z_s = z|\rho)]^2 \,>\,0$, which
is independent of $s>0$ because of
scaling/self-similarity properties of $(Z,\rho)$.
The scaling properties of $(Z,\rho)$ are also closely related to the
``aging'' of $(X,\tau)$.  Our main technical result is a general
convergence criterion for localization and aging functionals of 
diffusions/walks
$Y^{(\epsilon)}$ with (nonrandom) speed measures $\mu^{(\epsilon)} \to \mu$
(in a sufficiently strong sense).
\end{abstract}

\ep



\ep

%%%%%%%%%%%%%%%%%% SECTION 1 %%%%%%%%%%%%%%%%%%%%%%%%%%%%%%%%
\section{Introduction}
\setcounter{equation}{0}
\label{sec:int}

In this paper we continue the study of localization in the one-dimensional
Random Walk with Random Rates (RWRR), begun in \cite{kn:FIN}
(or equivalently of chaotic time dependence in the related Voter
Model with Random Rates (VMRR)---see below and \cite{kn:FIN}).
We also relate localization to ``aging'', a phenomenon of considerable
interest in out-of-equilibrium physical systems, such as glasses
(see, e.g., \cite{kn:BCKM} for a review).

\begin{defin} {\em (Random Walk with Random Rates)}
\label{def:rwrr}
The RWRR, $(X,\tau)$, is a continuous-time simple
symmetric random walk on $\z^d$, $X\, = \,(X_t : t\geq0,\,X_0=0)$,
where the time spent at site $i$ before
taking a step has an exponential distribution of mean $\tau_i$,
and where the $\tau_i$'s are i.i.d.~positive
random variables with common distribution $F$;
thus it is a random walk in the random
environment, $\tau = (\tau_i:i\in \z^d)$.
Except when otherwise noted, we restrict attention to $d=1$.
\end{defin}

When $F$ has a finite mean, it can be shown (e.g., by the convergence
results of \cite{kn:S}, as discussed below) that
(for a.e. $\tau$) there is a central limit theorem for $X_t$, and
more generally an invariance principle, i.e.,
that $\e X_{t/{\e^2}}$ converges to a Brownian motion
as $\e \to 0$. On the other hand, when $F$ has infinite mean
with a power law tail of exponent $\a <1$, one expects power law subdiffusive
behavior (with an exponent depending on both $\a$
and $d$); for reviews of the physics literature on subdiffusivity
in random environments, see, e.g., \cite{kn:HB,kn:Is}.
Logarithmic subdiffusivity \cite{kn:Sin}, as occurs in other
commonly studied random walks in random environments \cite{kn:Sol},
would presumably occur in an RWRR if the tail of $F$ were itself
logarithmic, but the more natural context for an RWRR is a power law
tail for $F$.

More striking than subdiffusivity,
and the main result of \cite{kn:FIN}, is that for $\a <1$ and
$d=1$,
there is localization in the sense that (for a.e. $\tau$) as $t \to \infty$,
\beq
\label{eq:loc1}
\sup_{i \in \z} \P(X_t = i|\tau) \not\to 0
\eeq
or equivalently
\beq
\label{eq:loc2}
\sum_{i \in \z} [\P(X_t = i|\tau)]^2 \not\to 0.
\eeq
An essential purpose of this paper is to relate this localization
to an appropriate scaling limit of $X$, in which it turns out
that Brownian motion
is replaced by a singular diffusion $Z$ (in a random
environment) --- singular here meaning that the single time
distributions of $Z$ are discrete. We remark that there is also
localization in the random walks of \cite{kn:Sol,kn:Sin}, as shown by
Golosov~\cite{kn:Go}, but both the localization and scaling limits
are of a somewhat different character there (as one would expect
in cases of logarithmic subdiffusivity);
for results about aging in these types of random walks
with random environments, see \cite{kn:DMF,kn:DGZ}.

Kawazu and Kesten~\cite{kn:KK} treated
the similar problem of finding the scaling
limit of a random walk with i.i.d. random bond
rates $\lambda_i$ (for transitions from $i$ to $i+1$
and from $i+1$ to $i$). Their random walk is in fact also
related to the VMRR and hence to our RWRR, with
$\lambda_i = 1/(2\tau_i)$. The scaling limit of \cite{kn:KK}
(see also \cite{kn:Sch1,kn:Sch2}) for $\alpha < 1$,
obtained by a similar approach based on \cite{kn:S} as
the one used here, is also a diffusion, but one that is
nonsingular
in the sense that the single time distributions
are continuous.
Our analysis of the type of localization
exhibited in (\ref{eq:loc1})--(\ref{eq:loc2})
(i.e., at {\it individual\/} points)
requires a stronger type of convergence to the scaling limit
than was needed in \cite{kn:KK}, as we explain later;
this strengthened convergence is the main new technical result
of the paper.

A convenient quantity, with which to express the relation
between localization and the scaling limit,
is the ``amount of localization'' at time t, as measured by
\beq
\label{eq:loc3}
  q_t = \E\sum_{i \in \z} [\P(X_t = i|\tau)]^2,
\eeq
where the expectation is with respect to $\tau$.
A main result of this paper, Theorem~\ref{teo:ordpar},
is that as $t \to \infty$,
$q_t$ converges to a (nonrandom) $q_\infty \in (0,1)$ (depending on
$\a  \in (0,1)$), which can itself be
expressed by a formula (see (\ref{eq:diffparam}) 
and (\ref{eq:ordpar}) below) analogous to
(\ref{eq:loc3})
with the singular diffusion $Z$ replacing the random walk $X$.

Our analysis of the scaling limit of $(X,\tau)$ will also yield results
about aging of the RWRR. As in the extensive physics literature on
the subject (see, e.g., \cite{kn:BCKM} and the references therein;
see \cite{kn:BDG} for rigorous work on aging in certain mean field
models),
we will consider a quantity
$R(t_w +t,t_w)$ that measures the behavior
of the system at a time $t_w +t$, after it has been aged for time
$t_w$. Normal aging corresponds to there being a well-defined nontrivial
limit function when $t$ and $t_w$ are scaled proportionally:
\beq
\label{eq:aging1}
{\cal R}(\theta)\,=\, \lim_{t_w\to\infty;t/t_w \to \theta} R(t_w +t,t_w).
\eeq
One interesting example of an $R$ for which such a limit follows from our
results is $R(t_w +t,t)=q_t (t_w)$, where
\beq
\label{eq:aging2}
q_t (t_w)\,=\, \E\sum_{i \in \z} [\P(X_{t_w +t} = i|\tau,X_{t_w})]^2.
\eeq
Of course, $q_t (0)=q_t$, corresponding to the amount of localization
after time $t$, starting from a fresh ($t_w=0$) system
with $X_0 = 0$ that has not been
aged. As with $q_\infty$, the limit function ${\cal R}(\theta)$ will be
given by a formula (see (\ref{eq:aging3})) like (\ref{eq:aging2}), but
with $X$ replaced by the diffusion $Z$. It follows from (\ref{eq:aging3})
that ${\cal R}(\theta)$ tends to $1$ as $\theta \to 0$ and to $q_\infty$
as $\theta \to \infty$.
Other examples of RWRR quantities that exhibit normal aging are
the (unconditional) probabilities
$\P(X_{t_w + t} =X_{t_w})$, which we discuss below, and
\beq
\label{eq:aging2.5}
\P(\max_{t_w \leq t' \leq t_w + t} \tau_{X_{t'}} >
  \max_{0 \leq t' \leq t_w} \tau_{X_{t'}}),
\eeq
which measures the prospects for ``novelty'' in this aging system.

Before explaining
more about $Z$ and it's random environment, we make a short digression
to point out that $q_\infty $ is a natural object of study
also for the related VMRR (as it is for other similar spin systems
with stochastic dynamics).

The one-dimensional (linear) VMRR is the continuous-time Markov process
$\st$ with state space $\{\s(i):i\in\z\}=\{-1,+1\}^{\z}$
in which, at rate $1/\tau_i$, site $i$ chooses (with equal probability)
one of it's two neighbors
(say $i'$) and replaces $\s(i)$ with $\s(i')$. The initial state
$\s_0$ is taken to be $\xi=(\xi_i:i\in\z)$, with the $\xi_i$'s i.i.d.~and
equally likely to be $+1$ or $-1$. Chaotic Time Dependence (CTD) is said
to occur if (conditional on $(\xi,\tau)$) the distribution of $\st$
has multiple subsequence limits as $t \to \infty$. (For a discussion of the
possible occurence of CTD in other more physical spin systems, see
\cite{kn:FIN, kn:NS1}.) Since the alternative to
CTD for this VMRR would be for the distribution to converge to the symmetric
mixture of the degenerate measures on the constant (identically $+1$ or
identically
$-1$) states, CTD is equivalent to the existence of some predictability
about the state for some arbitrarily large times, based on complete knowledge
of the inital state (and the environment of rates). In \cite{kn:FIN}, CTD
is proved to occur for a fat-tailed F (with $\a <1$) by showing that
(for a.e.~$(\xi,\tau)$ and every $k$) $\E[\s_t(k)|\xi,\tau]$ does not
converge as $t \to \infty$,
whereas the absence of CTD would require convergence to zero. A natural
quantity
measuring the amount of CTD/predictability
(see, e.g., \cite{kn:NNS}) is thus
\beq
\label{eq:ergavg}
\lim_{t \to \infty}\lim_{L \to \infty}(2L+1)^{-1}\sum_{k=-L}^{k=L}\E^2
[\s_t(k)|\xi,\tau] =
\lim_{t \to \infty}\E\{\E^2 [\s_t(0)|\xi,\tau]\}.
\eeq
But by the standard fact that a time-reversed voter model corresponds to
coalescing random walks, it easily follows, by doing the outermost expectation
first over $\xi$ and then over $\tau$, that
\beq
\label{eq:ordparam}
\E\{\E^2 [\s_t(0)|\xi,\tau]\} =
\E(\E\left\{\left.\E^2[\s_t(0)|\tau,\xi]\right|\tau\right\})
= \E\sum_{i \in \z} [\P(X_t = i|\tau)]^2 = q_t.
\eeq
Thus, in the VMRR,
the natural dynamical order parameter for CTD is just $q_\infty$.

Of course, it should be noted, that the existence of the
$t \to \infty$ limit in (\ref{eq:ergavg}) is not at all obvious---especially
in view of CTD.
(The $L \to \infty$ limit is a consequence of the spatial ergodicity of
$(\xi,\tau)$.) Indeed, we prove its existence 
(see Theorem~\ref{teo:ordpar}) by  expressing the $t \to
\infty$
limit of (\ref{eq:ordparam}) in terms of a scaling limit of $(X,\tau)$,
i.e., by  showing that as $t \to \infty$,
\beq
\label{eq:diffparam}
q_t  \to \E\sum_{z \in \r} [\P(Z_s = z|\rho)]^2 > 0,
\eeq
where $(Z,\rho)$ is a (singular) one-dimensional diffusion $Z$ in a random
environment
$\rho$. Here $s>0$ is arbitrary, and by the singularity of $Z$, we mean that
(conditional on $\rho$) the distribution of $Z_s$ is discrete, even though $Z$
is a bona-fide diffusion with continuous sample paths. We shall see why the
above
expression for $q_\infty$, which describes the amount of localization of
$(Z,\rho)$
at time $s$ does not in fact depend on $s$ (as long as $s \neq 0$),
a fact that may at first seem surprising (since $Z_s \to 0$ as
$s \to 0$, almost surely). Indeed this lack of dependence follows from
the scaling/self-similarity properties of $(Z,\rho)$ which imply
that (conditioned on $\rho$) the distribution of
$s^{\alpha/(\alpha + 1)}Z_s$ is a {\it random\/} measure on $\r$
whose distribution (arising from its dependence on $\rho$)
does not depend on $s>0$.
We now give a precise definition of this diffusion in a random
environment, $(Z,\rho)$.

\begin{defin} {\em (Diffusion with random speed measure, $(Z,\rho)$)}
\label{def:drsm}
The random environment $\rho$, the spatial scaling limit
of
the original environment $\tau$ of rates on $\z$, is a random discrete
measure,
$\sum_i W_i \delta_{Y_i}$, where the countable collection of $(Y_i,W_i)$'s
yields an
inhomogeneous Poisson point process on $\r \times (0,\infty)$ with density
measure
$dy\, \a w^{-1-\a} dw$. 
Conditional on $\rho$, $Z_s$ is a diffusion process
(with $Z_0 = 0$) that can be expressed as a time change of a standard 
one-dimensional Brownian
motion $B(t)$ with speed measure $\rho$, as follows~\cite{kn:IM}.
Letting $\ell(t,x)$ denote the local time at $x$ of $B(t)$, define
\beq
\label{eq:loctimint}
\phi_t^\rho:=\int\ell(t,y)\,d\rho(y)
\eeq
and the stopping time
$\psi_s^\rho$ as the first time $t$ when $\phi_t^\rho = s$ (so that
$\psi^\rho$
is the inverse function of $\phi^\rho$); then $Z_s = B(\psi_s^\rho)$. 
\end{defin}

Note that although $\rho$ is discrete,
the set of $Y_i$'s is a.s.~dense in $\r$ because the density measure is
non-integrable at $w=0$. For
(a deterministic) $s>0$,
the distribution of $Z_s$ is a discrete measure whose atoms are precisely
those of
$\rho$; this is essentially because the set of times when $Z$ is anywhere
else than
these atoms has zero Lebesgue measure.

The next theorem gives the limit (\ref{eq:diffparam})
as part of the convergence of the rescaled 
random walk $(X,\tau)$ to
the diffusion $(Z,\rho)$. The proof is not presented here because a more
complete result explaining the nature of this convergence is
provided later in Theorem \ref{teo:sca}. 

\bteo
\label{teo:ordpar}
Assume that $\P(\tau_0 > 0)=1$ and $\P(\tau_0>t)=L(t)/t^\a$, where
$L$ is a nonvanishing slowly varying function at infinity
and $\a<1$. Then for $\e > 0$, there exists $ c_\e > 0$ with
$  c_\e \to 0$ as $\e \to 0$, so that for any fixed $s>0$,
the distribution of $\Ze_s=\e X_{s/(c_\e\e)}$, conditioned on $\tau$
and thus regarded as a random probability measure on $\r$,
converges to the distribution of $Z_s$, conditioned on $\rho$,
in such a way that 
\beq
\label{eq:ordpar}
q_{s/(c_\e\e)} = \E\sum_{i\in\z} [P(\Ze_s=\e i | \tau)]^2 
\to \E\sum_{z\in\r} [P(Z_s=z | \rho)]^2.
\eeq
\eteo

We now return to a discussion of aging in the RWRR.
Analogously to (\ref{eq:diffparam}), we have ${\cal R}(\theta)$
of (\ref{eq:aging1})--(\ref{eq:aging2}) given by
\beq
\label{eq:aging3}
\lim_{t'\to\infty} q_{\theta t'} (t')\,=\,
\E\sum_{z \in \r} [\P(Z_{s+\theta s} = z|\rho,Z_s)]^2.
\eeq
The validity of this limit also follows from the results
and techniques of Sections 2 and 3
of the paper --- see Remark~\ref{rmk:aging}.
Here, the self-similarity properties of $(Z,\rho)$ imply that the RHS
of (\ref{eq:aging3}) depends only on $\theta$ and {\it not\/} on $s$
(for $0<s<\infty$), explaining the basic signature of normal aging ---
that the asymptotics of $q_t(t_w)$ depend only on the asymptotic
ratio of $t/t_w$.

Another example of an RWRR localization quantity with normal aging behavior is
\beq
\label{eq:aging4}
q'_t(t_w)\,=\, \P(X_{t_w + t} =X_{t_w})\,=\,
\E\, \P(X_{t_w + t} =X_{t_w}|\tau, X_{t_w}).
\eeq
In this case, the asymptotic aging function, ${\cal R}'(\theta)$,
would have limits of $1$ and $0$ respectively as $\theta \to 0$ and $\infty$.
Interestingly, a related quantity,
\beq
\label{eq:aging5}
q^*_t(t_w)\,=\, \P(X_{t_w + t'} =X_{t_w}\, \forall \, t' \in [0,t]),
\eeq
exhibits what is known as ``subaging'' (see, e.g., \cite{kn:RMB},
where a one-parameter family of models extending the RWRR are
studied nonrigorously, for general $d$). I.e.
(assuming, for simplicity, that the tail of $F$ satisfies
$u^\a\P(\tau_0 > u) \to K \in (0,\infty)$), there is a nontrivial limit
when $t/(t_w)^\eta \to \theta$ as $t_w \to \infty$, for some
$0<\eta<1$; here $\eta = 1/(1+\a)$ (for $0<\a<1$).
The difference in behavior between $q'$ and $q^*$ is due to the fact that
during the time interval $[t_w, t_w + \theta t_w]$, each visit
of the random walk to $X_{t_w}$ takes an amount of time of order
$t_w^{1/(1+\alpha)}$, but there are of order $t_w^{\alpha/(1+\alpha)}$
visits.
A related fact, in the scaling limit, is that for
$s,s'>0$, the diffusion process $Z$ has (for a.e.~$\rho$)
\beq
\label{eq:aging6}
\P(Z_{s+s'}=Z_s \,|\, \rho)\,>\,0 \quad \mbox{but} \quad
\P(Z_{s+s''}=Z_s \, \forall \, s''\in[0,s']\,\,|\, \rho)\,=\,0.
\eeq
This existence of different scaling regimes for {\it different\/}
quantities in a single model may be compared and contrasted to
the search for multiple scaling regimes in the {\it same\/} quantity
(see, e.g., \cite{kn:RMB}), where $R(t_w +\theta (t_w)^\eta ,t_w)$
and $R(t_w +\theta' (t_w)^{\eta'} ,t_w)$ with $\eta \neq \eta'$ would
both have nontrivial limits. (In fact, something weaker than this
is claimed in \cite{kn:RMB} for the $q^*_t(t_w)$ of (\ref{eq:aging5}).)

To see the lack of dependence of the RHS's of (\ref{eq:diffparam})
and (\ref{eq:aging3}) on $s$, we
may proceed
as follows. For $\lambda > 0$, consider the rescaled Brownian motion and
environment,
\beq
\label{eq:rescaling}
B^\lambda (t) = \lambda^{-1/2} B(\lambda t); \quad \rho^\lambda =
    \sum_i (\lambda^{-1/2})^{1/\a}W_i \delta_{\lambda^{-1/2} Y_i}.
\eeq
Since $B^\lambda$ and $\rho^\lambda$ are equidistributed with $B$ and $\rho$,
it follows that if we define a diffusion $Z^\lambda$ as the time-changed
$B^\lambda$
using speed measure $\rho^\lambda$, then $(Z^\lambda,\rho^\lambda)$ is
equidistributed
with the original diffusion in a random environment $(Z,\rho)$.
On the other hand,
on the original probability space on which $B$ and $\rho$ are defined,
one has $Z^\lambda_s = \lambda^{-1/2}Z_{\lambda^{(\a+1)/(2\a)}s}$, so that
the RHS's of (\ref{eq:diffparam}) and (\ref{eq:aging3})
remain the same when $s$ is replaced by
$\lambda^{(\a+1)/(2\a)}s$, and thus cannot depend on $s$.

To best understand how $(Z,\rho)$ arises as the scaling limit
of $(X,\tau)$, one should use
the fact that not only diffusions, but also random walks (or more accurately,
birth-death
processes) can be expressed as time-changed Brownian
motions~\cite{kn:S,kn:IM}.
In particular, if for any $\e>0$, we take as speed measure
\beq
\label{eq:discenv}
{{\rho}}^{(\e)}:=\sum_{i\in \z}c_\e\tau_i\d_{\e i},
\eeq
where the parameter $c_\e >0$ is yet to be determined, and then
do the time-change on the rescaled Brownian motion $B^{1/{\e^2}}$,
the resulting process is a rescaling of the original random walk $X$,
namely ${ {Z}}^{(\e)}_s = \e X_{s/(c_\e \e)}$. When the distribution
$F$ of the $\tau_i$'s has a finite mean, then by the Law of Large Numbers,
taking $c_\e = \e$, ${ {\rho}}^{(\e)}$ converges to (the mean of $F$ times)
Lebesgue measure and ${ {Z}}^{(\e)}$ converges to a Brownian motion
as $\e \to 0$~\cite{kn:S,kn:KV}. On the other hand, if $1-F(u) = L(u)/u^{\a}$
with $\a < 1$
and $L(u)$ is slowly varying at infinity~\cite{kn:F}, then by choosing 
$c_\e$
appropriately (as $\e^{1/\a}$ times a slowly varying function
at zero --- see (\ref{eq:coup4}) below) one has (from the classical
theories of domains of attraction and extreme value statistics) convergence
(in various senses, to be discussed) of ${ {\rho}}^{(\e)}$ to
the random measure $\rho$.

The idea that there should also follow
some kind of convergence of $({ {Z}}^{(\e)},{  {\rho}}^{(\e)})$ to
$(Z,\rho)$ should by now be quite clear. And indeed the basic convergence
results of~\cite{kn:S} are enough to imply, for example, that a
functional like
\beq
\label{eq:pseduoparam}
\E\{[P(a \leq { {Z}}_s^{(\e)} \leq b\,|\,{ {\rho}}^{(\e)})]^2\}
\eeq
(for {\it deterministic\/} $a,b$) converges to the corresponding quantity
for $(Z,\rho)$. But they are not sufficient to get localization
quantities like $q_{s/(c_\e \e)} =
\E\sum_{z \in \r} [P({ {Z}}_s^{(\e)} = z\, |\,{ {\rho}}^{(\e)})]^2$
to converge.
As mentioned earlier, the work of \cite{kn:KK} was also based on the
time-changed Brownian motion approach of \cite{kn:S,kn:IM},
but for their random walk and scaling limit, the convergence
results of \cite{kn:S} are sufficient.

The problem in our case
is not primarily with the randomness of ${ {\rho}}^{(\e)}$ (i.e., of $\tau$)
and $\rho$, but occurs already when
considering the nature of convergence of a process
$Y^{(\e)}(t)$ that is a Brownian motion time-changed with a {\it
deterministic\/}
speed measure $\mu^{(\e)}$. The convergence results of~\cite{kn:S} imply that
if $\mu^{(\e)} \to \mu$ vaguely, then (for example) one has weak convergence
of the distribution $\mbe$ of $Y^{(\e)}(t_0)$ to the corresponding
$\mb$. But we need stronger convergence.

This stronger convergence is the
subject of Section 2, which contains the main technical result
of the paper, Theorem~\ref{teo:conv},
in which weak convergence is combined with ``point process convergence''.
By point process convergence for (say) a discrete measure
$\sum_i \we_i \d_{y^{(\e)}_i}$ to $\sum_i w_i \d_{y_i}$
(where we have expressed each sum so that the atoms are not repeated),
we mean that the subset of $\r \times (0,\infty)$ consisting of all the
$(y^{(\e)}_i,\we_i)$'s converges to the set of all
$(y_i,w_i)$'s --- in
the sense that every open disk (whose closure is a compact
subset of $\r \times (0,\infty)$) containing exactly $m$ of the $(y_i,w_i)$'s
($m=0,1,\dots$) with none on its boundary,
contains also exactly $m$ of the $(y^{(\e)}_i,\we_i)$'s for all small $\e$.
Our
technical result is that vague plus point process convergence for the
speed measures $\mu^{(\e)} \to \mu$ implies the same for the distributions
at a fixed time $t_0$; i.e., $\mbe \to \mb$.

Going from this result for
a sequence of deterministic speed measures to our context of random speed
measures
requires a bit more work, which is presented in 
Sections \ref{sec:rwrr} and \ref{sec:scalimit} of the paper.
The way we handle that, which may be of independent interest, is to replace
the random measures ${{\rho}}^{(\e)}$ which only converge (in our two senses)
in distribution, by a different (but also natural) coupling
for the various $\e$'s than that
provided by the space of the original $\tau_i$'s so that convergence becomes
almost sure. This coupling is presented in Section~\ref{sec:rwrr} and
its convergence properties are given in Proposition~\ref{prop:coup}.
We note that almost sure convergence was also obtained in the
scaling limit results of \cite{kn:KK} by means of a coupling argument,
but there the coupling was an abstract one. In our situation, because
of the need to handle point process convergence, a concrete coupling
seems more suitable, in addition to being more natural.

We close the introduction by noting that we have restricted attention to the
scaling
limit of a single RWRR. In the context of the VMRR, which originally led to
our interest
in localization, one should consider the scaling limit of coalescing RWRR's.
Furthermore,
one should also study the scaling limit of the VMRR directly. These issues
will be taken
up in future papers.

%%%%%%%%%%%%%%%%%%%%%%%%%%%%%%%%%%%%%%%%%%%%%%%%%%%%%%%%%%%%%
%%%%%%%%%%%%%%%%%% SECTION 2 %%%%%%%%%%%%%%%%%%%%%%%%%%%%%%%%
%%%%%%%%%%%%%%%%%%%%%%%%%%%%%%%%%%%%%%%%%%%%%%%%%%%%%%%%%%%%%

\section{The continuity theorem}
\setcounter{equation}{0}
\label{sec:thm}

%%%%%%%%%%%%%% definitions of processes %%%%%%%%%%%%%%%%%%%%%%%%%%%%

Let $\mu,\,\me,\,\e\!>\!0,$ be non-identically-zero,
locally finite measures on $\r$.
%and let $\mu_d,\,\me_d,\,\e>0,$
%be their discrete parts. 
Let $Y_t,\,\Xe_t,\,t\!\geq\!0$,
$\Xe_0=Y_0=x,$
be the Markov processes in one dimension obtained by time
changing a standard Brownian motion
through $\mu,\me$; i.e., let
$B = B(s),\,s\geq0,$ be a standard Brownian
motion (with $B(0) = 0$) and let
\beq
\label{eq:tc1}
\phi_s(x):=\int\ell(s,y-x)\,d\mu(y),\quad
\psi(x)=\psi_t(x):=\phi_t^{-1}(x),\quad
Y_t=B(\psi_t(x))+x;
\eeq
\beq
\label{eq:tc2}
\fes(x):=\int\ell(s,y-x)\,d\me(y),\quad
\psi^{(\e)}(x)=\pset(x):=(\phi^{(\e)}_t)^{-1}(x),\quad
\xet=B(\pset(x))+x,
\eeq
where $\ell$ is the Brownian local
time of $B$~\cite{kn:S,kn:IM}.
Notice that,
since $\ell(s,y)$ is nondecreasing in $s$ for all $y$,
$\phi_s(x)$ and $\fes(x)$ are nondecreasing in $s$ and so their
(right-continuous)
inverses, $\psi_t(x)$ and $\pset(x)$, respectively, are well-defined.
Processes described in this way are known in the
literature as {\em quasidiffusions},
{\em gap diffusions} or
{\em generalized diffusions}~(\cite{kn:Ku,kn:Kn,kn:KW} and
references
therein). They generalize the usual diffusions in that the
{\em speed measures} $\mu$ can be zero in intervals, thus including
birth and death and other processes.

One fact about those processes we will need below is the following
formula from p.~641 of~\cite{kn:S}. Let $Y_0=x$;
for any Borel set $A$ of the reals,
\beq
\label{eq:form}
\int_0^t1\{Y_s\in A\}\,ds=\int_A\ell_Y(t,x,y)\,d\mu(y)
\eeq
almost surely, where $\ell_Y(t,x,y)=\ell(\psi_t(x),y-x)$.

%%%%%%%%%%%%%% definitions of types of convergence %%%%%%%%%%%%%%%%%%%%
We discuss now the types of convergence we will need for our results.
Let $\cm$ be the space of locally finite measures on $\r$ and $\cp$ 
its subspace of probability measures.
\begin{defin} {\em (Vague convergence)}
\label{def:vague}
Given a family $\nu,\,\nue,\,\e>0,$
in $\cm$,
%with their discrete parts
%$\nu_d,\,\nue_d,\,\e>0$, 
we say that $\nue$ converges {\em vaguely}
to $\nu$, and
write
$\nue\va\nu$, as $\e\to0$, if for all
continuous real-valued
functions $f$ on $\r$ with bounded support
$\int f(y)\,d\nue(y)\to\int f(y)\,d\nu(y)$ as $\e\to0$. 
\end{defin}
\begin{defin} {\em (Point process convergence)}
\label{def:pp}
For the same family, we say that 
$\nue$ converges {\em in the point process sense}
to $\nu$, and write
%$\nue_d\ppa\nu_d$ or 
$\nue\ppa\nu$, as $\e\to0$,
provided the following is valid:
if the atoms 
of $\nu$, $\nue$ are, respectively, at the {\it distinct\/} locations
$y_i$, $\xe_{i'}$
with weights $w_i$, $\we_{i'}$, then
the subsets $\ve\equiv\cup_{i'}\{(\xe_{i'},\we_{i'})\}$ of $\r\times(0,\infty)$
converge to
$V=\cup_i\{(y_i,w_i)\}$ as $\e\to 0$ in the sense that for any open
$U$ %\subset\r\times(0,\infty)$
whose closure
${\bar U}$ is a compact subset of $\r\times(0,\infty)$
such that its boundary contains no points of $V$,
the number of points $|\ve\cap U|$ in $\ve\cap U$
(necessarily finite since
$U$ is bounded and at a finite distance from $\r\times\{0\}$)
equals $|V\cap U|$ for all $\e$ small enough.
\end{defin}

%In the case that the family in Definition~\ref{def:vague}
%is in $\cp$, t
These notions can be related to the following
condition, where for $\nu \in \cp$ we 
order the $(y_i,w_i)$'s (the locations and weights of the 
atoms of $\nu$) so that $w_{i_1}\geq w_{i_2}\geq\ldots$, where 
$w_{i_1}$ is the largest weight, $w_{i_2}$ is the second largest, 
and so forth.
For a measure not in $\cp$, we use an arbitrary ordering of the atoms.
\bcond
\label{cond}
For each $l\geq1$, there exists $j_l(\e)$ such that
\beq
\label{eq:cond}
(y_{j_l(\e)},w_{j_l(\e)})\to(y_{i_l},w_{i_l})
\quad\mbox{as $\e\to0$}.
\eeq
\econd

We now establish a useful relationship among the above notions.

\bprop
\label{prop:rel}
For any family $\nu,\,\nue,\,\e>0,$ in $\cm$, the 
following two assertions hold.
If $\nue\ppa\nu$ as $\e\to0$, then Condition 1 holds.
If Condition 1 holds and $\nue\va\nu$ as $\e\to0$, 
then $\nue\ppa\nu$ as $\e\to0$.
\eprop
\noindent{\bf Proof.}

The first assertion is straightforward. Suppose the second one 
is false. According to the definitions above, that
means that there exists an open set $U_0$ whose closure is in
$\r\times(0,\infty)$ and a sequence $(\e_n)$ tending to $0$ as
$n\to\infty$
such that $|V^{(\e_n)}\cap U_0|\ne|V\cap U_0|$ for all $n$.
By Condition~\ref{cond} it must then be that
$|V^{(\e_n)}\cap U_0|>|V\cap U_0|$ for all large enough $n$.
That means that either
there exist ${\hat\imath}$, $w^\ast>0$ and
sequences  $(\e'_j)$, $(k_j)$ and $(k_j')$, with
$\e'_j\to0$ as $j\to\infty$ and $k_j\ne k_j'$ for all $j$,
such that
$y^{(\e'_j)}_{k_j},y^{(\e'_j)}_{k_j'}\to y_{\hat\imath}$,
$w^{(\e'_j)}_{k_j}\to w_{\hat\imath}$ and
$w^{(\e'_j)}_{k_j'}\to w^\ast$ as $j\to\infty$
or there exist a point
$(y^\ast,w^\ast)\in\r\times(0,\infty)\backslash V$
and sequences $(\e'_j)$ and $(k_j)$, with
$\e'_j\to0$ as $j\to\infty$,
such that $(y^{(\e'_j)}_{k_j},w^{(\e'_j)}_{k_j})\to(y^\ast,w^\ast)$
as $j\to\infty$.
In either case we get a contradiction to vague convergence of
$\nue$ to $\nu$ by taking a continuous
function $\tilde f$ that approximates sufficiently well the indicator
function of either $y_{\hat\imath}$ or $y^\ast$, depending on the case,
and showing that $\int\tilde f\,d\nu^{(\e_j)}$
is bounded below away from $\int\tilde f\,d\nu$.
\square

\ep

We leave it to the reader to find an example where Condition 1
holds but point process convergence does not. The following is
a useful corollary of Proposition~\ref{prop:rel}.

\bprop
\label{prop:sumsq}
Let $\nu,\,\nue,\,\e>0,$ be any family in $\cp$.  
If as $\e\to0$ both $\nue\ppa\nu$ and $\nue\va\nu$, 
then as $\e\to0$ 
\beq
\label{eq:sumsq}
\sum_{i'}[w^{(\e)}_{i'}]^2\to\sum_i[w_i]^2.
\eeq
\eprop
\noindent{\bf Proof.}

By the first assertion of Proposition~\ref{prop:rel},
Condition~\ref{cond} holds. This in turn implies that
\beq
\liminf_{\e\to0}\sum_j[w^{(\e)}_j]^2\geq\sup_k\sum_{l=1}^k[w_{i_l}]^2
=\sum_{i}[w_{i}]^2.
\eeq
%Also, $\nue\ppa\nu$ implies that
This together with the distinctness of the
$(y_i,w_i)$'s also implies that for any $k$
the indices $j_1(\e),\dots,j_k(\e)$ are distinct for small
enough $\e$. Furthermore, it implies that if $k$ and $\delta$
are such that $w_{i_k} > \delta > w_{i_{k+1}}$,
then for small enough $\e$
\beq
\sup\{w^{(\e)}_j:\,j \notin \{j_1(\e),\dots,j_k(\e)\}\} < \delta.
\eeq
To see this, note that otherwise along some subsequence
$\e = \e_l \to 0$ there would be an index $j^*(\e) \notin
\{j_1(\e),\dots,j_k(\e)\}$
with $\liminf y^{(\e)}_{j^*(\e)} \geq \delta$ and either
(i) $y^{(\e)}_{j^*(\e)} \to y^* \in (-\infty,+\infty)$ or else
(ii) $|y^{(\e)}_{j^*(\e)}| \to \infty$. Case (i) would
contradict $\nue\ppa\nu$, while case (ii) would imply that
the family $\{\nue\}$ is not tight, which would contradict
$\nue\va\nu$ since $\nue$ and $\nu$ are all probability measures.
Using the above choice of $k$ and $\delta$, we thus have
\beq
\limsup_{\e\to0}\sum_j[w^{(\e)}_j]^2\leq\sum_{l=1}^k[w_{i_l}]^2
+\limsup_{\e\to0}\sum_{j}\delta w^{(\e)}_j
=\sum_{l=1}^k[w_{i_l}]^2+\delta.
\eeq
Letting $k\to\infty$ and $\delta \to 0$ completes the proof. \square

\ep

We are ready to state the main result of this section; its proof
will begin after two corollaries are presented.

\bteo
\label{teo:conv}
Let $\me,\mu, Y^{(\e)}, Y$  be as above and fix any deterministic
$t_0>0$ and $x \in \r$.
Let $\mbe$  denote the distribution of
$\Xe_{t_0}$ (with $\Xe_{0} = x)$
%, let $\mbe_d$ denote its discrete part,
and define $\mbz$
%, \mbz_d$ 
similarly for $Y_{t_0}$.
Note that $\mbe = D_{t_0,x}(\me)$ and $\mbz = D_{t_0,x}(\mu)$, where
$D_{t_0,x}$ is some deterministic function from the non-identically-zero
measures in $\cm$ to $\cp$. Suppose
\beq
\label{eq:c1}
\me\va\mu\quad\mbox{and}\quad\me\ppa\mu\quad\mbox{as $\e\to0$}.
\eeq
Then, as $\e\to0$,
\beqn
\label{eq:c2}
\mbe\va\mbz\quad\mbox{and}\quad
\mbe\ppa\mbz.
\eeqn
\eteo

\brm
\label{rmk:aging}
To study limits involving two (or more) times
(see, e.g., (\ref{eq:aging2}), (\ref{eq:aging3}), (\ref{eq:aging4})),
some straightforward extensions of Theorem~\ref{teo:conv} are
useful. One of these is that (\ref{eq:c2}) remains valid if
$\Xe_{0} = x^{(\e)}$ with $x^{(\e)} \to x$. Another is that the
single-time distribution $\mbe$ of $\Xe_{t_0}$ can be replaced
by the multi-time distribution of $(\Xe_{t_1},\dots,\Xe_{t_m})$,
with point process convergence for measures on $\r^m$ defined in the
obvious way.
\erm

The following is an immediate consequence of Theorem~\ref{teo:conv}
and Proposition~\ref{prop:sumsq}.

\bcor
\label{cor:op}
Under the same hypotheses,
the weights of the atoms of $\mbe$ and $\mbz$ satisfy
\beq
\label{eq:op}
\sum_{j}[\bw^{(\e)}_j]^2
\ar\sum_{i}[\bw_i]^2
\quad\mbox{as $\e\to0$}.
\eeq
\ecor

\brm
\label{rmk:op}
More explicitly,~(\ref{eq:op}) takes the form
\beq
\label{eq:op1}
\sum_{y\in\r}[\P(\Ye_{t_0}=y)]^2
\ar\sum_{y\in\r}[\P(Y_{t_0}=y)]^2
\quad\mbox{as $\e\to0$},
\eeq
or, equivalently, if $Y^{(\e)\prime}_{t}$ (resp.~$Y^{\prime}_{t}$)
is an independent copy of $Y^{(\e)}_{t}$ (resp.~$Y_{t}$), then
\beq
\label{eq:op2}
\P(Y^{(\e)\prime}_{t_0}=\Ye_{t_0})
\ar\P(Y^{\prime}_{t_0}=Y_{t_0})
\quad\mbox{as $\e\to0$}.
\eeq
\erm

\ep

%%%%%%%%%%%%%% proofs  %%%%%%%%%%%%%%%%%%%%%%%%%%%%

Let
us also note that the above arguments yield another corollary. To state
it, we denote by $\cal{D}(\nu)$, for $\nu$ a probability
measure on $\r$, the $\{0,1,2,\dots,\infty\}$-valued measure on $(0,1]$
with ${\cal{D}}(\nu)(\Gamma)$ the number of
$x$'s in $\r$ such that $\nu(\{x\})\!\in\!\Gamma$; i.e., $\cal{D}(\nu)$
describes the set of all weights $w_i$ of the atoms of $\nu$, counting
multiplicity. Of course, since $\nu$ is a probability measure,
${\cal{D}}((\delta,1])<\infty$ for any $\delta>0$. The above arguments
show that $\mbe\va\mbz$ and $\mbz\ppa\mbz$ together imply
that ${\cal{D}} (\mbe)\Rightarrow{\cal{D}} (\mbz)$, where this latter 
convergence means that $\int f\,d{\cal D}(\mbe)\to\int f\,d{\cal D}(\mbz)$ for
any bounded continuous $f$ that vanishes in a neighborhood of the origin. 
Thus we have

\bcor
\label{cor:op2}
Under the same hypotheses, ${\cal{D}} (\mbe)\Rightarrow{\cal{D}} (\mbz)$
as $\epsilon\to0$.
\ecor

\ep

\noindent{\bf Proof of Theorem~\ref{teo:conv}.}

%%%%%%%%%%%%%% vague convergence %%%%%%%%%%%%%%%%%%%%%%%%%%%%

%\label{stone}
The vague convergence assertion in~(\ref{eq:c2}) is contained in
Corollary 1 of~\cite{kn:S}. Actually, the latter result is stronger.
It states that $\{\Ye_t,\,t\in[0,T]\}$ converges in distribution to
$\{Y_t,\,t\in[0,T]\}$ as a process (in the Skorohod topology),
$T>0$ arbitrary.
We will indeed use the stronger result in the argument for point
process convergence later on.
The fixed $t_0$ case is a rather simple and straightforward consequence
of the Brownian representation~(\ref{eq:tc1})-(\ref{eq:tc2}), so we
next briefly indicate an argument.

Since
$\ell(s,y)$ can be taken continuous in $(s,y)$
and of bounded support
in $y$ for each $s$, the first assumption in~(\ref{eq:c2})
implies that $\fes(x)\to\phi_{s}(x)$ as $\e\to0$ for all $s$. It follows
that  $\pset(x)\to\psi_{t}(x)$ as $\e\to0$ for all $t$ where
$\psi(x)$
is continuous. It suffices now to argue
that for any deterministic $t$, $\psi(x)$ is almost surely
continuous at $t$. For that, notice that $\psi(x)$ is
discontinuous at $t$ (if and) only if $\phi(x)$ has a
{\em plateau}
at height $t$, i.e., only if $\phi_{T+s}(x)-\phi_{T}(x)=0$ for some 
$s>0$, where
$T=\inf\{s'\geq0:\phi_{s'}(x)=t\}$. But, from the definition of
$\phi(x)$
and
monotonicity of $\ell$, that means that
\beq
\label{eq:abs}
\ell(T+s,y-x)-\ell(T,y-x)=0\quad\mbox{for $\mu$-almost every $y$}.
\eeq
Now, the definition of $T$ implies that
$\phi_{T-s'}(x)<t$ for all $s'>0$. This implies that $B(T)=y_0-x$ for
some $y_0$ in the support of $\mu$. But given that, since $T$ is a
stopping time, $\ell(T+s,y_0-x)-\ell(T,y_0-x)$ is distributed
like  $\ell(s,0)$ and thus is strictly positive for all
$s>0$. The continuity of $\ell$ now implies that there exists $\d>0$
such that $\ell(T+s,y-x)-\ell(T,y-x)>0$ if $|y-y_0|<\d$, which
contradicts~(\ref{eq:abs}). This settles the vague convergence assertion
in~(\ref{eq:c2}).

To prove the point process convergence of~(\ref{eq:c2}), 
by the second assertion of
Proposition~\ref{prop:rel} and the vague convergence 
of~(\ref{eq:c2}) just proven, 
it is enough to show that Condition~\ref{cond} holds. 
For that, we will need the following three lemmas.

\blem
\label{lem:loc1}
The set of locations of the atoms of $\mb$,
$\cup_i\{\by_i\}$, is contained in that of $\mu$, $\cup_{i'}\{y_{i'}\}$.
\elem

\noindent{\bf Proof.}

It is a result from the general theory of quasidiffusions~\cite{kn:Ku,kn:Kn}
that for a process $Y'$ living on a finite interval $I$
(with appropriate boundary conditions),
there exists
a symmetric continuous transition density $p'_I(t,x,y)$ which is
strictly positive and such that
\beq
\label{eq:trans}
\P(Y'_t\in\,dy|Y'_0=x)=p'_I(t,x,y)\,\mu(dy)\quad\mbox{for $t>0$, $x,y\in I$}.
\eeq
This would imply the result immediately if our process Y were such
a finite
interval process, but it is not.
However, if we condition on its history being contained
within a fixed interval, then we can use~(\ref{eq:trans}). The details are as
follows.

Let $A_{t,l}=\{Y_s\in [-l,l]\mbox{ for }s\leq t\}$, where $l>|x|$, $t\geq0$.
Then, on $A_{t,l}$, $\{Y_s,s\leq t\}$ is distributed like
$\{Y'_s,s\leq t\}$ on the analogous $A'_{t,l}$, where $Y'$ is the diffusion
with speed measure
$\mu':=\mu|_{(-l-1,l+1)}$ (and boundary conditions at $\pm(l+1)$ as
in~\cite{kn:Ku}). More precisely, there is a coupling between $Y, Y'$
and a third process $Y''$ with speed measure $\mu'$ and {\it killing\/}
boundary conditions at $\pm(l+1)$, such that $\{Y_s,s\leq t\}=\{Y'_s,s\leq t\}$
on $A''_{t,l}$. Thus
\beq
\label{eq:Loc1}
\P(Y_t=y_0|Y_0=x)- \e_{t,l}=p'_I(t,x,y_0)\,\mu'(y_0)-\e'_{t,l},
\eeq
where $I=[-l-1,l+1]$ and
$0 \leq \e_{t,l}, \e'_{t,l} \leq\P((A''_{t,l})^c)$.
If $\mu(y_0)=0$, then $\P(Y_t=y_0|Y_0=x)\leq\P((A''_{t,l})^c)$ for all
$l$. Then $\P(Y_t=y_0|Y_0=x)\leq\lim_{l\to\infty}\P((A''_{t,l})^c)=0$.
To obtain the vanishing of the last limit, we
first notice that, for given $T>0$,
$\P(A''_{t,l})\geq\P(B_s+x\in [-l,l],\,s\leq T)-\P(\psi_t(x)>T)$,
where $B$ is the standard Brownian motion in~(\ref{eq:tc1}).
The latter probability is bounded above by $\P(\phi_T(x)\leq t)$.
Thus $\liminf_{l\to\infty}\P(A''_{t,l})\geq
1-\limsup_{T\to\infty}\P(\phi_T(x)\leq t)$.
{}From~(\ref{eq:tc1}) and the known fact that almost surely
$\lim_{T\to\infty}\ell(T,x')=\infty$ for all $x'$,
the latter lim sup is seen to
vanish. The proof of the lemma is complete. \square

\blem
\label{lem:cont}
For all $y_0$, $\P(Y_t=y_0|Y_0=x)$ is continuous in $t>0$.
\elem

\noindent{\bf Proof.}

In view of Lemma~\ref{lem:loc1}, it suffices to consider the case where
$\mu(y_0)>0$. Let $t',t$ be such that $|t'-t|\leq1$. Imitating the argument
of the proof of that lemma,
\beqn
\nonumber
&|\P(Y_{t'}=y_0|Y_0=x)-\P(Y_{t}=y_0|Y_0=x)|&\\
&\leq|p'_I(t,x,y_0)-p'_I(t',x,y_0)|\,
\mu(y_0)\,+\,\P((A''_{t,l})^c) + \P((A''_{t',l})^c).&
\eeqn
Then $\lim_{t'\to t}|\P(Y_{t'}=y_0|Y_0=x)-\P(Y_{t}=y_0|Y_0=x)|
\leq\P((A''_{t,l})^c)+
\P((A''_{t+1,l})^c)$ for all $l$ and the result follows as in the proof of
Lemma~\ref{lem:loc1}. \square

\blem
\label{lem:loc2}
The set of locations of the atoms of $\mb$,
$\cup_i\{\by_i\}$, contains that of $\mu$, $\cup_{i'}\{y_{i'}\}$.
\elem

\noindent{\bf Proof.}
This is a corollary to the continuity lemma just given
and formula~(\ref{eq:form}). From
that formula, we have
\beq
\label{eq:f1}
\int_t^{t'}\P(Y_s=y_0|Y_0=x)\,ds=\E(\ell_Y(t',x,y_0)-
\ell_Y(t,x,y_0))\mu(y_0).
\eeq
We claim that the above expectation is strictly positive for all $x$, $y_0$ and
$t'>t$, if $\mu(y_0)>0$. This is a consequence of the definition
(see below (\ref{eq:form})) of $\ell_Y$ and the fact that there is
strictly positive probability that between
the two stopping times, $\psi_t(x)$ and $\psi_{t'}(x)$, the Brownian
motion $B$ will pass through $y_0 - x$ and hence will strictly
increase its local time there. Thus the integral on the left
hand side of~(\ref{eq:f1}) is also strictly positive.
This implies that for all $x$,
in every open interval of the
positive reals, there exists an $s$ such that
$\P(Y_s=y_0|Y_0=x)>0$. This and the continuity in time of these probabilities
imply that, in every interval $(0,t)$ there exists an $s$ such that
$\P(Y_{t-s}=y_0|Y_0=y_0)\P(Y_s=y_0|Y_0=x)>0$. By the Markov property
and time homogeneity of $Y$, the latter product
is a lower bound for $\P(Y_t=y_0|Y_0=x)$. The lemma follows. \square

\ep

We return to the proof of Condition~\ref{cond}.
Let $i_l$ be such that $\by_{i_l}=y_{i_l'}.$
Since $\me\ppa\mu$, there exists $j_l'(\e)$ so that
$(\ye_{j_l'(\e)},\we_{j_l'(\e)})\to(y_{i_l'},w_{i_l'})$ as $\e\to0$.
Since, by Lemmas \ref{lem:loc1} and \ref{lem:loc2},
$\{\by^{(\e)}_j\}=\{y^{(\e)}_{j'}\}$, we can define
$j_l(\e)$ so that $\by^{(\e)}_{j_l(\e)}=y^{(\e)}_{j'_l(\e)}$.
We already then have $\by^{(\e)}_{j_l(\e)}\to\by_{i_l}\,\,(=y_{i_l'}).$
To obtain Condition~\ref{cond}, we must show that also
$\bw^{(\e)}_{j_l(\e)}\to\bw_{i_l}$, i.e., that
$\P(\Ye_{t_0}=\by^{(\e)}_{j_l(\e)})\to\P(Y_{t_0}=\by_{i_l})$.

Let us simplify the notation a bit by setting $t_0=1$,
$\by^{(\e)}_{j_l(\e)}=y^{(\e)}_{j'_l(\e)}=y_\e$,
$\by_{i_l}=y_{i'_l}=y_0$, $w^{(\e)}_{j'_l(\e)}=w_\e$,
$w_{i'_l}=w_0$, $\bw^{(\e)}_{j_l(\e)}=\bw_\e$ and
$\bw_{i_l}=\bw_0$.
Thus, as $\e\to0$, we have $\me\va\mu$, $y_\e\to y_0$,
$w_\e\equiv\me(y_\e)\to\mu(y_0)=w_0$, and  we must show
that  $\bw_\e\equiv\P(\Ye_1=y_\e)\to\P(Y_1=y_0)=\bw_0$.

We also already know that $\Ye_1\to Y_1$ in distribution
(i.e., $\mbe\va\mb$). It follows that
$\limsup_{\e\to0}\P(\Ye_1=y_\e)\leq\P(Y_1=y_0)$, since otherwise
$\mbe\va\mb$ would be violated. So we only need to prove
\beq
\label{eq:co1}
\liminf_{\e\to0}\P(\Ye_1=y_\e)\geq\P(Y_1=y_0).
\eeq
By convergence in distribution,
\beqnn
\P(Y_1=y_0)\=\lim_{\d\to0}\P(y_0-\d<Y_1<y_0+\d)\\
  \le\lim_{\d\to0}\liminf_{\e\to0}\P(y_0-\d\leq\Ye_1\leq y_0+\d)\\
  \le\lim_{\d\to0}\left[\liminf_{\e\to0}\P(\Ye_1=y_\e)\left]
     \left[\limsup_{\e\to0}\,\fr{\P(y_0-\d\leq\Ye_1\leq y_0+\d)}
         {\P(\Ye_1=y_\e)}\right].\right.\right.
\eeqnn
Hence to prove~(\ref{eq:co1}), it suffices to show that
\beq
\lim_{\d\to0}\limsup_{\e\to0}\,
  \fr{\P(y_0-\d\leq\Ye_1\leq y_0+\d)}{\P(\Ye_1=y_\e)}\leq1
\eeq
or, equivalently, that
\beq
\label{eq:co2}
\lim_{\d\to0}\liminf_{\e\to0}\,
  \fr{\P(\Ye_1=y_\e)}{\P(y_0-\d\leq\Ye_1\leq y_0+\d)}\geq1.
\eeq

%%%%%%%%%%%%%%% (Ie)-(IVe)   %%%%%%%%%%%%%%%%%%%%%%%%%%%%

Given any small $\d>0$, we want to find a small
$\d'=\d'(\d)>\d$ with
$\d'\to0$ as $\d\to0$ and small $\T=\T(\d)$ with $0<\T<1$, such that
the following will be valid.
\bei
%%%%%%%%%% (Ie)
\item[($\mathrm {I}_\e$)]\label{ie}\,\,
$\displaystyle\lim_{\d\to0}\liminf_{\e\to0}\,
\P\left(y_0-\d'\leq\Ye_{1-\T}\leq y_0+\d'\left|
         y_0-\d\leq\Ye_1\leq y_0+\d\right)\right.=1$
%%%%%%%%%% (IIe)
\item[($\mathrm {II}_\e$)]\label{iie}\,\,
$\displaystyle\lim_{\d\to0}\liminf_{\e\to0}\,
\P\left(y_0-\d\leq\Ye_{1-\T}\leq y_0+\d\left|
         y_0-\d'\leq\Ye_{1-\T}\leq y_0+\d'\right)\right.=1$
%%%%%%%%%% (IIIe)
\item[($\mathrm {III}_\e$)]\label{iiie}\,\,
$\displaystyle\lim_{\d\to0}\liminf_{\e\to0}\,
\P\left(\inf_{t\in[1-\T,1]}\Ye_t \geq y_0-\d',\,
\sup_{t\in[1-\T,1]}\Ye_t\leq y_0+\d'
\left|\Ye_{1-\T}\in[y_0-\d,y_0+\d]\right)\right.=1$
%%%%%%%%%% (IVe)
\item[($\mathrm {IV}_\e$)]\label{ive}\,\,
$\displaystyle\lim_{\d\to0}\liminf_{\e\to0}\,
\P\left(\Ye_t=y_\e\mbox{ for some }t\in[1-\T,1]\left|
  \Ye_{1-\T}\in[y_0-\d,y_0+\d]\right)\right.=1$
%%%%%%%%%%
\eei

If ($\mathrm {I}_\e$)-($\mathrm {IV}_\e$) hold,
then~(\ref{eq:co2}) would be a consequence of
\beq
\label{eq:co3}
\lim_{\d'\to0}\limsup_{\e\to0}\sup_{s>0}\,
\P\left(\Ye_{s}\neq y_\e \mbox{ and }
\Ye_t\in[y_0-\d',y_0+\d']\mbox{ for all }t\in[0,s]\left|\Ye_{0}=y_\e 
\right)\right.=0.
\eeq
But this follows from the assumptions $\me\va\mu$ and $\me\ppa\mu$
(so that $\me(y_\e)\to\mu(y_{0})$) as $\e\to0$ and the following lemma.

\blem
\label{email}
For any open interval $I$ containing $\ye$,
\beq
\label{eq:email}
\P\left(\Ye_{s}\neq y_\e \mbox{ and }
\Ye_t\in I\mbox{ for all }t\in[0,s] \left|\Ye_{0}=y_\e \right)\right.\leq
1 - \fr{\me(y_\e)}{\me(I)}.
\eeq
\elem

\noindent{\bf Proof.}

The first step of the proof is similar to part of the proof of
Lemma~\ref{lem:loc1}, except here we use a coupling between the original
process $\Ye_t$ and a different process on the finite interval, namely
the process ${\tilde Y}_t$, whose speed measure is the finite measure
$\me \cdot 1_{I}$. (Basically, ${\tilde Y}_t$ has reflecting
boundary conditions at both endpoints of $I$.) Let
$A_{s,I}$ denote the event that $\Ye_t \in I$ for all $t \in [0,s]$;
then we take a coupling in which $\{\Ye_t, 0 \leq t \leq s\}$
equals $\{{\tilde Y}_t, 0 \leq t \leq s\}$ on $A_{s,I}$. Then the 
probability in
(\ref{eq:email}) equals
\beq
\label{eq:email2}
\P(\tilde{Y}_s \neq y_\e \mbox{ and } A_{s,I}|\tilde{Y}_0 = y_\e)
  \leq \P(\tilde{Y}_s \neq y_\e |\tilde{Y}_0 = y_\e)
  = 1 - \P(\tilde{Y}_s = y_\e |\tilde{Y}_0 = y_\e).
\eeq
The proof is completed by applying the following lemma with $Y$ replaced by
$\tilde{Y}$, $\mu$ by $\me \cdot 1_{I}$, and $y$ by $y_\e$. \square

\blem
\label{lem:spectral}
Let $s \geq 0$ and $y \in \r$; then
\beq
\label{eq:email3}
\P(Y_s = y |Y_0 = y) \geq  \fr{\mu(y)}{\mu(\r)}.
\eeq
\elem

\noindent{\bf Proof.}

We may assume that $\mu(y)>0$ and $\mu(\r) < \infty$, since otherwise
the claim is trivially true. 
To avoid technical considerations about generators of quasidiffusions and
their spectral properties, we will temporarily further assume that
$\mu$ is {\em finitely supported}. Then $Y$ is a Markov jump process
with {\em finite} state space (the atoms of $\mu$) and has $[\mu(\r)]^{-1}\mu$
as its unique invariant distribution. Let $T_t$ denote the transition
semigroup acting on the finite-dimensional space, $L^2 (\r,d\mu)$:
\beq
\label{eq:email4}
T_t\,:\, f(x) \longmapsto \E(f(Y_t)|Y_0 =x).
\eeq
Then $T_t = e^{t{\cal L}}$, where 
the generator ${\cal L}$ has a simple eigenvalue $0$, with normalized 
eigenvector
the constant function $\Phi(x)=[\mu(\r)]^{-1/2}$, and the rest of
its spectrum strictly negative. Let $\Psi(x)$ be the (normalized in 
$L^2 (\r,d\mu)$)
function $[\mu(y)]^{-1/2}1_y$. Then by the spectral theorem, 
and denoting by $(\cdot,\cdot)$ the inner product in $L^2 
(\r,d\mu)$,
\beq
\label{eq:email5}
\P(Y_s = y |Y_0 = y) = ([\mu(y)]^{-1}1_y,T_s 1_y)=
(\Psi,e^{t{\cal L}}\Psi) = |(\Psi,\Phi)|^2 + \int_{-\infty}^{0-} e^{sl}d\nu(l),
\eeq
where $\nu$ (the spectral measure of $\Psi$ restricted to $(-\infty,0)$) is
a finitely supported non-negative measure on $(-\infty,0)$. It follows
that $\P(Y_s = y |Y_0 = y)$ is non-increasing in $s$
and converges, as $s\to\infty$, to $|(\Psi,\Phi)|^2=\mu(y)/\mu(\r)$.

We have now proved~(\ref{eq:email3}) when $\mu$ is finitely supported.
If not, we take a sequence $\mu^{[n]}$ of finitely supported measures
such that $\mu^{[n]}\va\mu$ as $n\to\infty$, with
$\mu^{[n]}(y)=\mu(y)>0$ and $\mu^{[n]}(\r)=\mu(\r)<\infty$
for all $n$. The corresponding processes $Y^{[n]}$ converge to $Y$
in distribution as $n\to\infty$ by the results of~\cite{kn:S} so that 
by standard weak convergence arguments,
\beq
\label{eq:email6}
\P(Y_s=y|Y_0=y)\geq\limsup_{n\to\infty}\P(Y^{[n]}_s=y|Y^{[n]}_0=y).
\eeq
The proof is completed by using~(\ref{eq:email3}), as already proved for
$(Y,\mu)$ replaced by $(Y^{[n]},\mu^{[n]})$.  \square

\ep

%%%%%%%%%%%%%%%%%  (Ie)-(IVe) hold   %%%%%%%%%%%%%%%%%%%%%%%%%%%%

It remains to show that ($\mathrm {I}_\e$)-($\mathrm {IV}_\e$) hold
(for some $\d'(\d)$ and $\T(\d)$).  From
the convergence in distribution (in the Skorohod topology) of the
processes (\cite{kn:S}; see the first paragraph of this proof,
following Corollary~\ref{cor:op2}), we
have, for example, that
\beqn
\nonumber
&\liminf_{\e\to0}
\P\left(\Ye_{1-\T}\in[y_0-\d',y_0+\d']\right)
\geq\P\left(Y_{1-\T}\in(y_0-\d',y_0+\d')\right)&\\
&=\lim_{\d''\ua\d'}
\P\left(Y_{1-\T}\in[y_0-\d'',y_0+\d'']\right),&
\eeqn

\beq
\limsup_{\e\to0}
\P\left(\Ye_{1}\in[y_0-\d,y_0+\d]\right)
\leq\P\left(Y_{1}\in[y_0-\d,y_0+\d]\right),
\eeq

\beqn
\nonumber
&\liminf_{\e\to0}&\!\!
\P\left(\Ye_{t}\in[y_0-\d',y_0+\d']\mbox{ for all }t\in[1-\T,1]\right)\\
&\geq\lim_{\d''\ua\d'}&\!\!
\P\left(Y_{t}\in[y_0-\d'',y_0+\d'']\mbox{ for all }t\in[1-\T,1]\right),
\quad\mbox{etc}.
\eeqn
Thus, ($\mathrm {I}_\e$)-($\mathrm {III}_\e$) are seen to follow
from the corresponding ($\mathrm {I}$)-($\mathrm {III}$) with
$\Ye$ replaced by $Y$ (and $\liminf_{\e\to0}$
deleted). For ($\mathrm {IV}_\e$), a different argument is required,
because the usual notions of convergence in distribution of
$\Ye$ to $Y$ will not work directly for
($\mathrm {IV}_\e$) and the analogous ($\mathrm {IV}$).
Instead, we replace ($\mathrm {IV}_\e$) by a stronger condition
($\mathrm {IV}'_\e$), stated in terms of the Brownian motion $B$
of (\ref{eq:tc1})--(\ref{eq:tc2}):
\bei
%%%%%%%%%% (IV'e)
\item[($\mathrm {IV}'_\e$)]\label{iv'e}\,\,
$\displaystyle\lim_{\d\to0}\liminf_{\e\to0}\,
\P\left(Q_{\e,x,[1-\T,1]}\leq y_0-\d,\,
Q^{\e,x,[1-\T,1]} \geq y_0+\d
\left| \Ye_{1-\T}\in[y_0-\d,y_0+\d]\right.\right)=1$,
\eei
where
\beq
\label{eq:min}
Q_{\e,x,[a,b]}=\inf_{s\in[\pse_{a}(x),\pse_b(x)]}(B(s)+x)
\eeq
and $Q^{\e,x,[a,b]}$ is defined analogously with inf replaced by sup.
This condition is stronger because $\me(y_\e)>0$ and so $Y^{(\e)}$
cannot skip over $y_\e$.
Now, as above, by the convergence in distribution of
$\Ye$ to $Y$, it suffices to prove the
corresponding condition ($\mathrm {IV}'$) for $Y$.

\ep

%%%%%%%%%%%%%%%%%  (I)-(IV') hold   %%%%%%%%%%%%%%%%%%%%%%%%%%%%t

It remains to show that ($\mathrm {I}$)-($\mathrm {III}$) and
($\mathrm {IV'}$) hold for some $\d'(\d)$ and $\T(\d)$.

Since the distribution of $Y_{t_0}$ has an atom at $y_0$,
it follows that
\beq
\P(Y_{t_0}=y_0|Y_{t_0}\in[y_0-\d'',y_0+\d''])\to1\mbox{ as }
\d''\to0
\eeq
for each $t_0$. From this and Lemma~\ref{lem:cont}
(and the vague continuity in $t$ of the distribution of
$Y_t$, from, e.g., \cite{kn:S}), ($\mathrm {II}$) follows
(provided $\d'\to0$ as $\d\to0$).

Similarly, assuming $\d'\to0$ as $\d\to0$, we can replace
($\mathrm {I}$) by
\bei
%%%%%%%%%% (I')
\item[($\mathrm I'$)]\label{i'}\quad
$\displaystyle\lim_{\d\to0}\,
\P\left(Y_{1-\T}=y_0|Y_{1}=y_0\right)=1$.
%%%%%%%%%%
\eei
But this follows, assuming $\T\to0$ as $\d\to0$, again from
the continuity of
$\P(Y_{t}=y_0)$ in $t > 0$.

Let us take ($\mathrm {IV'}$) now. The probability there
is bounded from below by
\beq
\label{eq:iv'}
\inf_{x\in{\rm supp}\mu\cap[y_0-\d,y_0+\d]}
\P\left.\left(\inf_{s\in[0,\psi_{\T}(x)]}B(s)+x \leq y_0-\d,\,
\sup_{s\in[0,\psi_{\T}(x)]}B(s)+x \geq y_0+\d \right|\Ye_{0}=x\right).
\eeq
Let $S'_\d$ denote the time it takes
for $B(s)+y_0-\d$ to first reach $y_0+\d$ and then come back to
$y_0-\d$. Then the expression in~(\ref{eq:iv'}) is bounded from
below by
\beq
\label{eq:iv'1}
%\P(T'_\d<\T).
\P(\phi_{S'_\d}(y_0-\d) < \T).
\eeq
Now $S'_\d\to0$ as $\d\to0$ almost surely. From the almost sure
continuity
of $\ell$, we have $\ell(S'_\d,y-(y_0-\d))\to\ell(0,y-y_0)\equiv0$
as $\d\to0$
and from this and the monotonicity
of $\ell$ in $t$
and the fact that $\ell(t,\cdot)$ has compact support almost surely
for every $t$, it follows straightforwardly that
$\phi_{S'_\d}(y_0-\d)\to0$ as $\d\to0$ almost surely.
That means that the probability in~(\ref{eq:iv'1}) would
tend to $1$  as
$\d\to0$ for any fixed $\T>0$ (i.e., not depending on $\d$).
But then we can choose a sequence $\T=\T(\d)$, with $\T\to0$ as $\d\to0$,
such that it still tends to $1$ as $\d\to0$. This establishes
($\mathrm {IV'}$).

Finally we need to choose $\d'$ so that
($\mathrm {III}$) is valid. The argument is analogous to the above one
for ($\mathrm {IV'}$). The probability in ($\mathrm {III}$)
is bounded from below by
\beq
\label{eq:iii}
\inf_{x\in{\rm supp}\mu\cap[y_0-\d,y_0+\d]}
\P\left.\left(\inf_{t\in[0,\T]}Y_t \geq y_0-\d',\,
\sup_{t\in[0,\T]}Y_t\leq y_0+\d'\right|Y_{0}=x\right).
\eeq
Let $Y'_t$ be a copy of $Y_t$, starting at $y_0-\d$
at time $0$ and let $Y''_t$ be a copy of
$Y_t$ starting at
$y_0+\d$ at time $0$. Let $T'_{\d\d'}$ and $T''_{\d\d'}$
denote the time it takes
for $Y'_t$ and $Y''_t$ to first reach beyond $(y_0-\d',y_0+\d')$,
respectively, and let $\T = \T (\delta)$ be as chosen
in the previous paragraph with $\T \to 0$ as $\delta \to 0$.
Then the expression in~(\ref{eq:iii}) is bounded from
below by
\beq
\label{eq:iii1}
%\P(T'_{\d\d'}\wedge T''_{\d\d'}>\T)=
\P(T'_{\d\d'}>\T) + \P(T''_{\d\d'}>\T) - 1.
\eeq
Let us take the first term
of~(\ref{eq:iii1}).
Consider $S'_{\d\d'}$, the quantity corresponding to
$T'_{\d\d'}$ for $B(s)+y_0-\d$.
Then, if $\d'$ were fixed, we would have that
$S'_{\d\d'}\to S'_{0\d'}$ as $\d\to0$ almost surely.
Now $T'_{\d\d'}>\T$ if
$\phi_{S'_{\d\d'}}(y_0-\d) > \T$.
By the same reasoning as above, we have, for fixed $\d'$, that
$\phi_{S'_{\d\d'}}(y_0-\d)\to\phi_{S'_{0\d'}}(y_0)$ as $\d\to0$
almost surely.
That means that, as $\d\to0$, the first term
of~(\ref{eq:iii1}) is bounded below by
$\P(\phi_{S'_{0\d'}}(y_0)>0)=1$
for any fixed $\d'>0$, since $S'_{0\d'}>0$ for any $\d'>0$ and
$\phi_{s}(y_0)>0$ for all $s>0$ almost surely.
That means that  $\P(T'_{\d\d'}>\T(\delta))\to1$ as $\d\to0$
for any fixed $\d'>0$ and the same can be argued analogously
for the second term
of~(\ref{eq:iii1}).
Thus we can choose a sequence $\d'=\d'(\d)$, with $\d'\to0$ as $\d\to0$,
such that
(\ref{eq:iii1}) tends
to $1$ as $\d\to0$. This establishes
($\mathrm {III}$) and thus Condition 1. Theorem~\ref{teo:conv}
follows. \square

%%%%%%%%%%%%%%%%%%%%%%%%%%%%%%%%%%%%%%%%%%%%%%%%%%%%%%%%%%%%%
%%%%%%%%%%%%%%%%%% SECTION 3 %%%%%%%%%%%%%%%%%%%%%%%%%%%%%%%%
%%%%%%%%%%%%%%%%%%%%%%%%%%%%%%%%%%%%%%%%%%%%%%%%%%%%%%%%%%%%%

\section{A coupling for the scaled random rates}
\setcounter{equation}{0}
\label{sec:rwrr}

As discussed briefly at the end of Section~\ref{sec:int},
the rescaled random walk with random rates, 
${ {Z}}^{(\e)}_{\cdot} = \e X_{\cdot /(c_\e \e)}$
is a quasidiffusion whose (random) speed measure
${{\rho}}^{(\e)}$, given by (\ref{eq:discenv}), only
converges {\it in distribution\/} to the (random)
speed measure $\rho$ of the scaling limit diffusion $Z$.
To take advantage of the results of Section~\ref{sec:thm},
it is convenient to find random measures $\tre$ equidistributed 
(for each $\e$) with ${{\rho}}^{(\e)}$ and such that
$\tre$ converges {\it almost surely\/} as $\e \to 0$ to 
$\rho$, in both the vague and point process senses of 
Section~\ref{sec:thm}. This will be done in this section
by constructing $\tau^{(\e)}$, equidistributed with $\tau$
for each $\e > 0$, on the natural probability space where
$\rho$ is defined.

Consider the L\' evy process (see, e.g., \cite{kn:Re,kn:ST,kn:Sa})
$V_x,\,x\in\r,\,V_0=0$, with stationary
and independent increments given by
\beq
\label{eq:gf}
\E\left[e^{ir(V_{x+x_0}-V_{x_0})}\right]=
e^{\a x\int_0^\infty(e^{irw}-1)\,w^{-1-\a}\,dw}
\eeq
for any $x_0 \in \r$ and $x \geq 0$. It satisfies
\beq
\label{eq:ytail}
\lim_{y\to\infty} y^\alpha\P(V_1 > y) = 1
\eeq
(\cite{kn:F}, Theorem XVII.5.3).
Let $\rho$ be the (random) Lebesgue-Stieltjes measure on the
Borel sets of $\r$ associated to $V$, i.e.,
\beq
\label{eq:mu}
\rho((a,b])=V_b-V_a,\,a,b\in\r,\, a<b,
\eeq
where we have chosen the process $V$ to have sample paths
that are right-continuous (with left-limits). Then
\beq
\label{eq:dmu}
\frac{d\rho}{dx}=\frac{dV}{dx}=\sum_j\,w_j\,\d(x-x_j),
\eeq
where the (countable) sum is over the indices of an
inhomogeneous Poisson point process $\{(x_j,w_j)\}$
on $\r\times(0,\infty)$ with density $dx\,\a w^{-1-\a}\,dw.$

For each $\e>0$, we want to define, in the fixed probability space
on which $V$ and $\rho$ are defined,
a sequence $\te_i,\,i\in\z,$ of independent random variables such that
\beq
\label{eq:coup1}
\te_i\sim\tau_0\quad\mbox{for every $i\in\z$}
\eeq
(where $\sim$ denotes equidistribution)
and with the following property: for a given family of constants
$c_\e,\,\e>0$, let
\beq
\label{eq:coup2}
\tre:=\sum_{i= -\infty}^\infty c_\e\te_i\d_{\e i};
\eeq
we demand that constants $c_\e$ can be chosen so that
\beq
\label{eq:coup3}
\tre\va\rho\quad\mbox{and}\quad\tre\ppa\rho\quad\mbox{as $\e\to0$,
almost surely.}
\eeq

Before specifying our construction of $\te$ in general, we
consider the very special case  where $\tau_0$ is equidistributed
with the positive $\alpha$-stable random variable $V_1$.
Note then that according to (\ref{eq:ytail}),
$\P(\tau_0>t)=L(t)/t^\a$ with $L(t) \to 1$ as $t \to \infty$.
Here, we may simply choose $c_\e = \e^{1/\alpha}$ and
take $\te$ to be the sequence of scaled increments of the 
L\' evy process $V$:
\beq
\label{eq:special}
\te_i=\frac1c_\e\,\!\left(V_{\e(i+1)}-V_{\e i}\right).
\eeq
The validity of (\ref{eq:coup1}) and (\ref{eq:coup3}) are
then elementary exercises, which we leave to the reader.
If $\tau_0$ is not $\a$-stable but
$t^\a \P(\tau_0>t) \to K \,\in \,(0,\infty)$ as $t \to \infty$,
one can still take $c_\e$ proportional to $\e^{1/\alpha}$, but
a more complicated definition of $\te_i$ may be needed;
without such an assumption on the 
distribution of $\tau_0$, $c_\e$ may also require
a more complicated definition. 

The next proposition extends the construction of $\te$ to quite
general distributions of $\tau_0$ by choosing
the following $c_\e$ and $\te_i$'s:
\beqn
\label{eq:coup4}
c_\e\=\left(\inf\{t\geq0:\P(\tau_0>t)\leq \e\}\right)^{-1}\\
\label{eq:coup5}
\te_i\=\frac1c_\e\,g_\e\!\left(V_{\e(i+1)}-V_{\e i}\right),
\eeqn
where $g_\e$ is defined as follows. Let $G:[0,\infty)\to[0,\infty)$
satisfy
\beq
\label{eq:coup6}
\P(V_1>G(x))=\P(\tau_0>x)\quad\mbox{for all $x\geq 0$.}
\eeq
$G$ is well-defined since $V_1$ has a continuous distribution.
Notice that that $G$ is
nondecreasing and right-continuous and thus has a
nondecreasing and right-continuous
generalized inverse $G^{-1}$. Let $g_\e:[0,\infty)\to[0,\infty)$
be defined as
\beq
\label{eq:coup7}
g_\e(x)=c_\e G^{-1}\!\left(\e^{-1/\a}x\right)\quad\mbox{for all $x\geq 0$.}
\eeq

\bprop
\label{prop:coup}
Suppose that $\P(\tau_0 > 0)=1$ and $\P(\tau_0>t)=L(t)/t^\a$, where
$L$ is a nonvanishing slowly varying function at infinity
and $\a<1$. Then
(\ref{eq:coup1}) and~(\ref{eq:coup3}) hold for $c_\e$ and $\te_i$ as
in~(\ref{eq:coup4})-(\ref{eq:coup7}).
\eprop

\ep

\noindent{\bf Proof.}
%We will prove~(\ref{eq:coup1}) now and postpone~(\ref{eq:coup3}) until
%later in this section (following Theorems \ref{teo:sca} and \ref{teo:op1}).
%
To establish~(\ref{eq:coup1}), by the stationarity of the increments of $V$,
it suffices to take $i=0$. Then, for $\e>0$
\beq
\P(\te_0>t)=\P(g_\e(V_\e)>c_\e t)=\P(V_\e>g_\e^{-1}(c_\e t))=
\P(V_\e>\e^{1/\a}G(t)),
\eeq
where $g_\e^{-1}$ is the right-continuous inverse of $g_\e$,
and we have used the
easily checked fact that
$g_\e^{-1}(\cdot)=\e^{1/\a}G(\cdot/c_\e)$.
The desired result (\ref{eq:coup1}) now follows
by the scaling relation $V_\e\sim\e^{1/\a}V_1$
(see~(\ref{eq:gf})) and~(\ref{eq:coup6})). 

It remains to derive (\ref{eq:coup3}) to finish
the proof of Proposition~\ref{prop:coup}. For that we
will need two main lemmas, as follows.

\blem
\label{ml1}
For any fixed $y>0$, $\ge(y)\to y$ as $\e\to0$.
\elem

\blem
\label{ml2}
For any $\d'>0$, there exist constants $C'$ and $C''$ in $(0,\infty)$
such that
$$
g_\e(x)\leq C'x^{1-\d'}\mbox{ for }\e^{1/\a}\leq x\leq1
\mbox{ and }\e\leq C''.
$$
\elem

The proofs of these two main lemmas are based on the
following four subsidiary
lemmas, whose proofs are given later.

\ep

\blem
\label{A}
$\displaystyle\fr1\e\,\P\left(\tau_0>\fr1{c_\e}\right)\to 1$ as $\e\to0$.
\elem

\blem
\label{A'} For $y>0$,
$\displaystyle\fr1\e\,\P\left(\tau_0>\fr y{c_\e}\right)\to\fr{1}{y^\a}$
as $\e\to0$.
\elem

\blem
\label{B} For any $\l>0$,
$\displaystyle\fr{c_\e}{c_{\l\e}}\to \l^{-1/\a}$
as $\e\to0$ and thus (by standard results, as in~\cite{kn:F})
$c_\e=\e^{1/\a}\tilde{L}(\e^{-1})$, where $\tilde{L}$ is a
positive slowly varying function at infinity.
\elem

\blem
\label{C}
There exists $\l>0$ sufficiently small such that
$ G^{-1}(y)\leq1/{c_{\l/y^\a}}$ for $y\geq1$
or, equivalently,
$ g_\e(x)\leq{c_\e}/{c_{\l\e/x^\a}}$ for $x\geq\e^{1/\a}.$
\elem

\noindent{\bf Proof of Lemma \ref{ml1}}

Let $g_\e^{-1}$ be the right-continuous generalized inverse of $g_\e$.
To prove $g_\e(y)\to y$, it suffices to prove that $g_\e^{-1}(y)\to y$.

%Now $G^{-1}(V_1)\stackrel{\mathrm dist}{=}\tau_0$, so
Now $G^{-1}(V_1) \sim \tau_0$, so
$g_\e(\e^{1/\a}V_1)=c_\e G^{-1}(\e^{-1/\a}\e^{1/\a}V_1)
%\stackrel{\mathrm dist}{=}c_\e\tau_0$, so $\P(\tau_0>y/c_\e)=$
\sim c_\e\tau_0$, and thus $\P(\tau_0>y/c_\e)$ equals
\beq
\label{eq:ml11}
\P(c_\e\tau_0>y)=\P(g_\e(\e^{1/\a}V_1)>y)
=\P(\e^{1/\a}V_1>g_\e^{-1}(y))=\P(V_1>\e^{-1/\a}g_\e^{-1}(y)).
\eeq
By~(\ref{eq:ytail}),
\beq
\label{eq:ml12}
\e^{-1}\P(V_1>\e^{-1/\a}y)\to 1/y^\a
\eeq
as $\e\to0$.
By~(\ref{eq:ml11}) and Lemma~\ref{A'},
\beq
\e^{-1}\P(V_1>\e^{-1/\a}g_\e^{-1}(y))=\e^{-1}\P(\tau_0>y/c_\e)\to
1/y^\a
\eeq
as $\e\to0$. This implies that
$\P(V_1>\e^{-1/\a}g_\e^{-1}(y))/\P(V_1>\e^{-1/\a}y)\to 1$
as $\e\to0$ and this plus~(\ref{eq:ml12}) implies that
$\limsup_{\e\to0} g_\e^{-1}(y)\leq y$ and
$\liminf_{\e\to0} g_\e^{-1}(y)\geq y$, completing the proof of
Lemma~\ref{ml1}. \square

\ep

\noindent{\bf Proof of Lemma~\ref{ml2}}

By Lemmas \ref{B} and \ref{C}, for $x\geq\e^{1/\a}$,
\beq
\label{eq:ml21}
g_\e(x)\leq\l^{-1/\a}x\fr{\tilde{L}(\e^{-1})}{\tilde{L}((x^\a/\l)\e^{-1})}
\eeq
for $\l>0$ small enough; the value of $\lambda$ will be chosen later.
We now use a result from p.~274 of~\cite{kn:F} about slowly varying
functions, stating that
$\tilde{L}(x)=a(x)e^{\int_1^x\fr{\D(y)}y\,dy}$, where $a(x)\to c\in(0,\infty)$
as $x\to\infty$ and $\D(y)\to 0$ as $y\to\infty$. The quotient in the right
hand
side of~(\ref{eq:ml21}) then becomes
\beq
\label{eq:ml22}
\fr{a(\e^{-1})}{a((x^\a/\l)\e^{-1})}\,
\exp\left\{\int_{(x^\a/\l)\e^{-1}}^{\e^{-1}}\fr{\D(y)}y\,dy\right\}.
\eeq
If $\e \leq \l$ so that $(x^\a/\l)\e^{-1} \geq 1/\l \geq \e^{-1}$, then the
absolute value of the latter integral is bounded above by
\beq
\label{eq:ml23}
\d\left|\int_{(x^\a/\l)\e^{-1}}^{\e^{-1}}\fr1y\,dy\right|\leq
\d|\log(x^\a/\l)|,
\eeq
where $\d=\d(\l)=\sup\{|\D(y)|,y>1/\l\}$,
and  thus the exponential in~(\ref{eq:ml22}) is bounded above
(for $\l \leq 1$, $x \leq 1$) by
\beq
\label{eq:ml24}
\l^{-\d}x^{-\a\d}.
\eeq
Thus, given $\d'>0$, we choose $\l \in (0,1)$ such that $\a\d(\l)\leq\d'$
and such that
$a(y) \in [c/2,c]$ for $y\geq\l^{-1}$.
The lemma now follows from~(\ref{eq:ml21})-(\ref{eq:ml24}) with
$C'=4 \l^{-(1+\d')/\a}$ and $C'' = \l$. \,\square

\ep

To complete the proof of our two main lemmas, it
remains to prove the subsidiary Lemmas \ref{A}, \ref{A'}, \ref{B} and \ref{C}.

\ep

\noindent{\bf Proof of Lemma \ref{A}}

By the definition (\ref{eq:coup4}) of $c_\e$,
$\P(\tau_0>c_\e^{-1})\leq\e$ and $\P(\tau_0>x)>\e$ for all $x< c_\e^{-1}$.
Thus, if the statement of the lemma is not true,
then there must exist $\d \in (0,1)$
and a sequence $(\e_i)$ with $\e_i>0$ for all $i$ and $\e_i\to0$ as
$i\to\infty$ such that $\P(\tau_0>c_{\e_i}^{-1})\leq\d\e_i$ for all $i$. But
then, given $\d'$ such that $\d^{1/\a}<\d'<1$, we have that
$\P(\tau_0>\d'c_{\e_i}^{-1}) > \e_i$ and so
\beq
\label{eq:svar0}
\fr{\P(\tau_0>\d'c_{\e_i}^{-1})}{\P(\tau_0>c_{\e_i}^{-1})}\geq\d^{-1}
%\fr{$\P(\tau_0>\d'c_{\e_i}^{-1})}{$\P(\tau_0>c_{\e_i}^{-1})}\geq\d^{-1}
\eeq
for all $i$. Since $c_{\e_i}^{-1}\to\infty$ and
$\P(\tau_0>\cdot)$ is regularly varying at infinity
(with exponent $-\a$), it follows that for any $\l > 0$,
\beq
\label{eq:svar}
\lim_{t \to \infty} \fr{\P(\tau_0>\l t)}{\P(\tau_0> t)} = \l^{-\a},
\eeq
which contradicts (\ref{eq:svar0}) since
$(\d')^{\a}>\d$. \square

\ep

\noindent{\bf Proof of Lemma \ref{A'}}

This is a consequence of Lemma~\ref{A}, the fact that
$c_\e^{-1}\to\infty$
as $\e\to0$, and ({\ref{eq:svar}), from which it follows that
$$ \fr{\P(\tau_0>y/c_\e)}{\P(\tau_0>1/c_\e)}\to\fr1{y^\a}.$$

\ep

\noindent{\bf Proof of Lemma \ref{B}}

By Lemma~\ref{A}, $(\l\e)^{-1}\P(\tau_0>1/c_{\l\e})\to 1$ or
equivalently $\e^{-1}\P(\tau_0>1/c_{\l\e})\to \l $
as $\e\to0$ while, by Lemma~\ref{A'},
$\e^{-1}\P(\tau_0>y/c_{\e})\to 1/y^\a.$
This implies that taking $y^\a=\l^{-1}$ that
$c_\e\l^{1/\a}/c_{\l\e}\to 1$ or
$c_\e/c_{\l\e}\to \l^{-1/\a}$ as $\e\to0$.

\ep

\noindent{\bf Proof of Lemma \ref{C}}

To show that $G^{-1}(y)\leq z$, it is enough to show that
$G(z)>y$. Thus we want to prove that $G(1/c_{\l/y^\a})>y$
for $y\geq1$ and some $\l>0$. By the definition (\ref{eq:coup6})
of $G$, $G(x)>y$
would be a consequence of $\P(V_1>y)>\P(\tau_0>x)$, where we
take $x=1/c_{\l/y^\a}$. Now there exists $K>0$ such that
$\P(V_1>y)>K/y^\a$ for $y\geq1$ (by~(\ref{eq:ytail})),
so it suffices to show that
$\P(\tau_0>1/c_{\l/y^\a})\leq K/y^\a$ for $y\geq1$ and some
$\l>0$; or, equivalently, taking $\e=\l/y^\a$, it suffices to
show that for some $\l>0$ and all $\e\leq\l$,
$\P(\tau_0>1/c_{\e})\leq K\e/\l$, or
$\P(\tau_0>1/c_{\e})/\e\leq K/\l$. By Lemma~\ref{A}, we may choose
$\l$ small enough so that for $\e\leq\l$,
$\P(\tau_0>1/c_{\e})/\e\leq 2$ and also small enough that
$K/\l\geq2$.

\ep

\noindent{\bf Completion of Proof of Proposition~\ref{prop:coup}.}
We still have to prove~(\ref{eq:coup3}).
This will be done using our two main lemmas \ref{ml1} and \ref{ml2}.
The point process convergence of (\ref{eq:coup3}) would follow
straightforwardly if we knew that
$g_\e(x_\e)\to x_0$ as $\e\to0$ whenever
$x_\e\to x_0>0$. To obtain that, due to the monotonicity
and right continuity of $g_\e(\cdot)$, it suffices that
$g_\e(y)\to y$ as $\e\to0$ for any fixed $y>0$, and that is given by
Lemma~\ref{ml1}.

We argue next why the vague convergence of (\ref{eq:coup3})
follows by using both Lemma~\ref{ml1}
and Lemma~\ref{ml2}.
Let $f$ be a continuous function with bounded support $I$. Then
\beq
\label{eq:p1}
\int f\,d\tre=\sum_{i\in\e^{-1}I}f(\e i)g_\e(V_{\e(i+1)}- V_{\e i}).
\eeq
For $y>0$, let $J_y=\{i\in\e^{-1}I: V_{\e(i+1)}- V_{\e i}\geq y\}$.
To estimate~(\ref{eq:p1}), we treat separately the sums over
$J_\d$, $J_{\e^{1/a}}\backslash J_\d$ and $\e^{-1}I\backslash J_{\e^{1/a}}$,
with $\d>\e^{1/a}$. From
Lemma~\ref{ml1}, it follows that as $\e \to 0$,
\beq
\label{eq:p2}
\sum_{i\in J_\d}f(\e i)g_\e(V_{\e(i+1)}- V_{\e 
i})\to\sum_{j:w_j\geq\d}f(x_j)w_j ,
\eeq
where $\cup_j\{(x_j,w_j)\}$ is the Poisson point process of~(\ref{eq:dmu}).

By Lemma~\ref{ml2}, we have that, given $\d'>0$ small enough
(to be chosen shortly),
for some finite constant $C$,
\beq
\label{eq:p3}
\sum_{i\in J_{\e^{1/a}}\backslash J_\d}f(\e i)g_\e(V_{\e(i+1)}- V_{\e i})\leq
C\sum_{i\in J_{\e^{1/a}}\backslash J_\d}(V_{\e(i+1)}- V_{\e i})^{1-\d'}.
\eeq
The latter sum is bounded above by
\beq
\label{eq:p4}
W_\d := \sum_{j:x_j\in I,w_j\leq\d}w_j^{1-\d'}.
\eeq
With $\d'>0$ chosen small enough so that $\d'+\a<1$, we claim that
$W :=\lim_{\d\to0}W_\d=0$ almost surely. Indeed, note that $W$ is well defined
in any case by monotonicity and is of course non-negative. We also have,
by a standard Poisson process calculation, that
\beq
\label{eq:p5}
\E(W_\d)=|I|\int_0^\d w^{1-\d'}w^{-1-\a}\,dw<\infty
\eeq
for all $\d>0$ and $\E(W_\d) \to0$ as $\d\to0$.
By dominated convergence, $\E(W)=0$ and the claim follows.

Finally, by the definition (\ref{eq:coup7}) of $ g_\e $
and its monotonicity, we have that
$ g_\e(x)\leq g_\e(\e^{1/a})=C c_\e $ for $x\leq\e^{1/a}$, where $C$ is some
finite constant. It then follows that
\beq
\label{eq:p6}
\sum_{i\in\e^{-1}I\backslash J_{\e^{1/a}}}f(\e i)g_\e(V_{\e(i+1)}- V_{\e i})
\leq C'c_\e\sum_{i\in\e^{-1}I} 1\leq C''c_\e\e^{-1}\to0
\eeq
as $\e\to0$, by Lemma~\ref{B}, since $\a < 1$.

Combining the above estimates, we  get that $\int f\,d\tre$
converges to
\beq
\label{eq:p7}
\lim_{\e\to0}\sum_{i\in\e^{-1}I}f(\e i)g_\e(V_{\e(i+1)}- V_{\e i})
=\lim_{\d\to0}\sum_{j:w_j\geq\d}f(x_j)w_j=\sum_{j}f(x_j)w_j=\int f\,d\rho.
\,\, \square
\eeq

%%%%%%%%%%%%%%%%%%%%%%%%%%%%%%%%%%%%%%%%%%%%%%%%%%%%%%%%%%%%%
%%%%%%%%%%%%%%%%%% SECTION 4 %%%%%%%%%%%%%%%%%%%%%%%%%%%%%%%%
%%%%%%%%%%%%%%%%%%%%%%%%%%%%%%%%%%%%%%%%%%%%%%%%%%%%%%%%%%%%%

\section{Scaling limit for the random walk with random rates}
\setcounter{equation}{0}
\label{sec:scalimit}

Let $X_t,\,t\geq0,\,X_0=0$, be a continuous time random walk
on $\z$ with inhomogeneous rates given by $\l_i=\tau_i^{-1}$,
$i\in\z$, where $\tau_i,\,i\in\z$, are i.i.d.~random variables
such that $\P(\tau_0 > 0)=1$ and $\P(\tau_0>t)=L(t)/t^\a$, where
$L$ is a nonvanishing slowly varying function at infinity
and $\a<1$.

We consider now the scaling limit of the random walk $X_t$.
Let
\beq
\label{eq:sca}
\Ze_t=\e X_{t/(c_\e\e)},\quad t\geq0.
\eeq
To study the limit
of $\Ze$, in
the presence of the random rates, which themselves converge vaguely
and in the point process sense, but only in distribution,
we will need a weak notion of vague and point process convergence,
as follows. Let $\tc$ be the class of bounded real functions $f$ on 
the space $\cal P$ of probability measures on $\r$
that are weakly continuous in the
sense that $f(\mu_n)\to f(\mu)$ as $n\to\infty$ for all
$\mu,\,\mu_n$, $n\geq1$, in $\cal P$ such that {\it both\/} $\mu_n\va\mu$
{\it and\/} $\mu_n\ppa\mu$ as $n\to\infty$.

%\begin{defin}
%\label{def:dw}
%Let $P,\,P_\e,\,\e>0$ be probability measures on $\cal M$.
%We say that $P_\e$ converges {\em doubly weakly} to $P$, denoted
%$P_\e\dw P$, if, as $\e\to0$,
%\beq
%\label{eq:dw}
%\int f\,dP_\e\to\int f\,dP\,\mbox{ for all } f\in\tc.
%\eeq
%We also use the notation $\pi^{(\e)} \dw \pi$ for random measures
%$\pi^{(\e)}$ and $\pi$, whose distributions are $P_\e$ and $P$ respectively.
%\end{defin}

Let $Z_t$ be the (random) quasidiffusion $Y_t$ as
in~(\ref{eq:tc1})-(\ref{eq:tc2}) above, but with speed measure $\mu$
taken to be the 
(random) discrete measure $\rho$ of (\ref{eq:mu})--(\ref{eq:dmu})
associated with the L\'evy process $V$.
For $t_0>0$
fixed, let $\rb$ and $\rbe$ be the (random) probability distributions of
$Z_{t_0}$ and $\Ze_{t_0}$, respectively; i.e., $\rb$ is the
conditional distribution of $Z_{t_0}$ given $\rho$ while $\rbe$
is the conditional distribution of $\Ze_{t_0}$ given $\tau$.
We can now state the following theorem, which 
is a consequence of Proposition~\ref{prop:coup} and
Theorem~\ref{teo:conv}.

\bteo
\label{teo:sca}
As $\e\to0$,
\beq
\label{eq:sca1}
%\rbe\dw\rb.
\E(f(\rbe)) \to \E(f(\rb))
\eeq
for all $f\in\tc$; in particular, 
%\eteo
%
%\bteo
%\label{teo:op1}
\beq
\label{eq:op3}
\E\sum_{x\in\r}\left[\rbe(\{x\})\right]^2
=\E\sum_{x\in\r}[\P(\Ze_{t_0} = x |\tau)]^2
\to\E\sum_{x\in\r}\left[\rb(\{x\})\right]^2
=\E\sum_{x\in\r}[\P(Z_{t_0} = x |\tau)]^2.
\eeq
\eteo

\noindent{\bf Proof of Theorem~\ref{teo:sca}}

$\Ze$ is distributed as a standard Brownian motion time changed through
the speed measure $\re$ (see~(\ref{eq:discenv}) and the 
beginning of Section~\ref{sec:thm}).
Let $\rho$ and $\tre$ be as in~(\ref{eq:mu}) and~(\ref{eq:coup2})
and let
${\tilde{Z}^{(\e)}}$ be a standard Brownian motion time changed through
$\tre$. By  Proposition~\ref{prop:coup},
\beq
\label{eq:iddist}
(\Ze,\re) \sim {(\tilde{Z}^{(\e)}},\tre).
\eeq
To obtain (\ref{eq:sca1}), it is thus enough, by~(\ref{eq:iddist})
and dominated
convergence, to show that $\bar{\tilde{\rho}}^{(\e)}$, the probability
distribution of $\tilde{Z}^{(\e)}_{t_0}$ (which is 
$D_{t_0,0}({\tilde{\rho}}^{(\e)})$ in the
notation of Theorem~\ref{teo:conv} and is random because
of its dependence on $\tre$ and hence on the L\'evy process $V$),
satisfies: $\bar{\tilde{\rho}}^{(\e)}\va\rb$ and
$\bar{\tilde{\rho}}^{(\e)}\ppa\rb$ almost surely. But that follows from
Proposition~\ref{prop:coup} and Theorem~\ref{teo:conv}. 

Then (\ref{eq:op3}) follows from~(\ref{eq:sca1}) with the function
$f$ on $\cp$ defined by
$f(\mu)=\sum_{x\in\r}\left[\mu(\{x\})\right]^2$, which 
belongs to $\tc$ by Proposition~\ref{prop:sumsq}.
\square

%CHANGE THE PROOF OF THIS PART OF THM TO USE PROP FROM SEC. 2.
%Similarly, it suffices to prove that almost surely
%\beq
%\label{eq:op4}
%\sum_{x\in\r}\left[\bar{\tilde{\rho}}^{(\e)}(\{x\})\right]^2
%\sum_{x\in\r}\left[\bar{\tilde{\rho}}^{(\e}}(\{x\})\right]^2
%\to\sum_{x\in\r}\left[\rb(\{x\})\right]^2
%\eeq
%and that follows from Proposition~\ref{prop:coup} and
%Corollary~\ref{cor:op}. \square

\ep

%%%%%%%%%%%%%%%% ACKNOWLEDGEMENTS %%%%%%%%%%%%%%%%%%%%%%%%%%%%%%
\noindent{\bf Acknowledgements.}

The authors have benefitted greatly from discussions with
Anton Bovier about aging and with Raghu Varadhan
about quasidiffusions. They also thank
Maury Bramson, Amir Dembo, Alice Guionnet, Harry Kesten, Tom Kurtz,
Enzo Marinari, Henry McKean and Alain Sznitman for useful
comments and suggestions. In addition, they thank the referee
for constructive pointers that led to an improved presentation.

The research of L.F.~is part of FAPESP theme project no.~99/11962-9
and CNPq PRONEX project no.~41/96/0923/00. It was also supported in
part by CNPq research grant no.~300576/92-7, and travel grants from FAPESP
no.~00/02773-7 and CNPq/NSF no.~910116/95-4.

The research of M.I.~ was supported in part by MURST under grants
{\em Stime, ordinamenti stocastici e processi aleatori} and
{\em Processi stocastici con struttura spaziale}.

The research of C.M.N.~was supported in part by the N.S.F.~under
grants DMS-9803267 and DMS-0104278.

Each of the authors thanks the departments and universities of
their coauthors for hospitality and support provided during
various visits. 

\ep

%%%%%%%%%%%%%%%%%%%%%% REFERENCES %%%%%%%%%%%%%%%%%%%%%%%%%%%%%%%

%\input bib

\vskip 5truemm
\parindent -20pt
\leftline{Instituto de Matem\'atica e Estat\'\i stica ---
Universidade de S\~ao Paulo}
\leftline{Cx.\ Postal 66281 --- 05315-970 S\~ao Paulo SP --- Brasil}
\leftline{{\tt<lrenato@ime.usp.br>}}

\vskip 3truemm
\parindent -20pt
\leftline{Dipartimento di Interuniversitario Matematica ---
Universit\` a di Bari}
\leftline{Via E. Orabona 4 --- 70125 Bari --- Italia}
\leftline{{\tt<isopi@dm.uniba.it>}}

\vskip 3truemm
\parindent -20pt
\leftline{Courant Institute of Mathematical Sciences ---
New York University}
\leftline{251 Mercer Street --- New York, NY 10012 --- USA}
\leftline{{\tt<newman@cims.nyu.edu>}}

\end{document}